\documentclass[preprint]{elsarticle}
\usepackage{hyperref}
\usepackage[centertags]{amsmath}
\usepackage{amsfonts}
\usepackage{amssymb}
\usepackage{diagbox}
\usepackage{adjustbox}
\usepackage{amsmath, bm}
\usepackage[ruled]{algorithm2e}
\usepackage{enumitem}
\usepackage{multicol}
\usepackage{multirow}
\usepackage{makecell}
\usepackage{ragged2e}
\usepackage{changes}
\usepackage{color}
\usepackage{ulem}

\newtheorem{remark}{Remark}
\journal{Journal}
\numberwithin{table}{section}
\usepackage[ruled]{algorithm2e}
\SetKwBlock{GrO}{Grouping optimization:}{end}
\SetKwBlock{PeO}{Perform optimization steps:}{end}

\numberwithin{equation}{section}
\numberwithin{figure}{section}
\begin{document}
\begin{frontmatter}
\title{An extrapolation-driven network architecture for physics-informed deep learning}

\author[1,2,3]{Yong Wang} 
\author[1,3]{Yanzhong Yao\corref{cor1}}
\ead{yao\_yanzhong@iapcm.ac.cn}
\author[1,3]{Zhiming Gao}

\cortext[cor1]{Corresponding author}
\address[1]{Institute of Applied Physics and Computational
Mathematics, Beijing 100088, China}
\address[2]{Graduate School of China Academy of Engineering Physics, Beijing 100088, China}
\address[3]{National Key Laboratory of Computational Physics, Beijing 100088, China}

\begin{abstract}
Current \textcolor{black}{physics-informed neural network (PINN) implementations} with sequential learning strategies \textcolor{black}{often experience} some weaknesses, such as the failure to reproduce the previous training results \textcolor{black}{when using} a single network, the difficulty to strictly ensure continuity and smoothness at the time interval nodes \textcolor{black}{when using}  multiple networks, and the increase in complexity and computational overhead. To \textcolor{black}{overcome}  these shortcomings,
we first investigate the extrapolation capability of the PINN method for time-dependent PDEs. Taking advantage of this extrapolation property, we generalize the training result obtained in
\textcolor{black}{a specific} time subinterval to 
\textcolor{black}{larger intervals} 
by adding a correction term to the network parameters of the subinterval. The correction term is determined by further training with the sample points in the added subinterval. Secondly, by designing an extrapolation control function with special characteristics and combining it with \textcolor{black}{a} correction term, we construct a new neural network architecture whose network parameters are coupled with the time variable, which we call the extrapolation-driven network architecture. Based on this architecture, using a single neural network, we can obtain the overall PINN solution of the whole domain with the following two characteristics: (1) it completely inherits the local solution of the interval obtained from the previous training, (2) at the interval node, it strictly maintains the continuity and smoothness that the true solution has. The extrapolation-driven network architecture allows us to divide a large time domain into multiple subintervals and solve the time-dependent PDEs one by one in \textcolor{black}{a} chronological order. This training scheme respects the causality principle and effectively overcomes the difficulties of the conventional PINN method in solving the evolution equation on a large time domain. Numerical experiments verify the performance of our  method. 
The data and code accompanying this paper are available at \href{https://github.com/wangyong1301108/E-DNN}{https://github.com/wangyong1301108/E-DNN}.
\end{abstract}

\end{frontmatter}

\section{Introduction}
Investigating and understanding the dynamics of complex physical processes is critical to science and engineering. 
\textcolor{black}{The dynamic behavior of objects is typically formulated by time-dependent partial differential equations (PDEs).} Providing efficient and highly accurate numerical methods for time-dependent PDEs is a goal that has been relentlessly pursued by many scientific computing researchers. In the past decades, researchers have made great progress in exploring the dynamics by using the traditional numerical methods, such as the finite difference method (FDM), the finite volume method (FVM), and the finite element method (FEM).  
\textcolor{black}{However, there are still some challenges in using these classical analytical or computational tools.} Traditional approaches not only face the curse of dimensionality, but also struggle to provide universal and effective frameworks for solving evolution equations with complex features.

\textcolor{black}{Deep neural networks (DNNs) have garnered widespread attention as a tool that adapts to high-dimensional data and has an innate capacity to model nonlinear relationships.} 
Based on DNNs, researchers have developed many deep learning methods~\cite{2018eweinan,2018SIRIGNANO20181339,RAISSI1,WANG2024113112}.
Among these methods, those using physics-informed neural networks (PINNs)~\cite{RAISSI1,w712178} have emerged as a popular option due to their simplicity and effectiveness~\cite{kim2021dpm}, \textcolor{black}{and PINN methods have already been used to address various problems in different fields~\cite{KISSAS2020112623,XU2023108900,YAO2023116395,taufik2023neural,10zishiyingditi}}.
Despite these successes, \textcolor{black}{they are} still in the early stages \textcolor{black}{in tackling} realistic problems with PINNs~\cite{WANG202111}. 

\textcolor{black}{We focus here} on solving time-dependent PDEs by PINNs.
In Ref.~\cite{PENWARDEN2023112464}, Penwarden et al.~argue that the conventional PINN method faces some difficulties when \textcolor{black}{trained to predict} over the entire spatio-temporal domain.

They summarize these difficulties into 3 types: (1) Zero solution, where PINNs tend to converge to a zero solution when the zero solution could minimize residual loss; 
(2) No propagation, where PINNs fail to propagate any information due to an insufficient number of residual points, resulting in invalid predictions;(3) Incorrect propagation, where PINNs quickly converge to an incorrect\textcolor{black}{/trivial} solution, perhaps because the optimizer falls into a local minimum \textcolor{black}{corresponding to a trivial soluton}.
To address these difficulties, researchers have proposed many improved versions of the conventional PINN method, 
\textcolor{black}{including domain decomposition strategies~\cite{JAGTAP2020113028,CiCPXPINN,pfoasd1010638830},  sequential learning strategies~\cite{PENWARDEN2023112464,CiCP-29-930,MATTEY2022114474,GUO2023112258},  and causal weights~\cite{WANG2024116813}.
While the method in~\cite{WANG2024116813}  has been successful in simulating complex systems, its training time  is very long. As a positive development, Ref.~\cite{eaFDSFSDge} mitigates the problem by modifying the causality parameter. Additionally,} sequential learning strategies are often preferred in applications~\cite{Physicsinformedneuralnetwork,CHEN2024111423} because they naturally respect the \textit{principle of causality}~\cite{WANG2024116813} and have \textcolor{black}{a relatively higher} computational efficiency.

\textcolor{black}{
The extrapolation properties of  multilayer perceptrons (MLPs)
~\cite{158898,balestriero2021learninghighdimensionamounts} and PINNs~\cite{Karniadakis,HAGHIGHAT2021113741,9933886}  have been widely studied.
Prior work demonstrates that MLPs usually succeed in interpolation, but fail to extrapolate well~\cite{xu2021how,ZHU2023116064}.
However, as mentioned in Refs.~\cite{Karniadakis,HAGHIGHAT2021113741}, PINNs are capable of interpolation and extrapolation by using physical information as soft penalty constraints.
Ref.~\cite{9933886} illustrates that PINNs have credible interpolation and extrapolation properties when applied to velocity models containing gaps.
Designing appropriate neural network architectures based on the extrapolation property of PINNs is a feasible  way to address the above issues.}

For sequential learning strategies, there are two alternative network options: a single feed-forward fully connected neural network~\cite{PENWARDEN2023112464,CiCP-29-930,MATTEY2022114474,GUO2023112258} or multiple feed-forward fully connected neural networks~\cite{PENWARDEN2023112464,CiCP-29-930}. 
The main implementation steps of using a single network are as follows: (1) First train the neural network on a small time subinterval, which is called pre-training. It is well known that training on a small domain generally yields more accurate predictions~\cite{10191822}; (2) Using the resulting parameters of the small subinterval as initial guesses, then train the neural network on a larger interval including the pre-training subinterval. 
\textbf{The weakness of the single-network approach is that the PINN solution obtained on the large interval cannot accurately reproduce that of the small subinterval.}
\textcolor{black}{One option to address the weakness is by using neuron splitting~\cite{ewqrewrJB023703}.}
On the other hand, the main implementation steps of using multiple feed-forward fully connected neural networks are as follows:
(1) Divide the whole time domain into several time subintervals, each subinterval is expressed by a separate neural network;
(2) Use some prediction values of previous subintervals as supervised data, train the neural network of each subinterval sequentially one by one.
(3) Combine the prediction functions of all subintervals to form the PINN solution for the whole domain.
\textbf{The disadvantage of the multi-network approach is that, \textcolor{black}{it is difficult to strictly maintain the continuity and smoothness at the interval nodes when combining the PINN solutions of multiple neural networks as a potential solution of the evolution equations in the whole time domain.}
In addition, the implementation of such methods is very complex and computationally time-consuming~\cite{10.1007/978-3-031-08754-7_45}.}

This paper attempts to address the above issues from the perspective of updating the neural network architecture. 
By analyzing the extrapolation property of the conventional PINN method, we construct a kind of extrapolation-driven neural network architecture to compensate for the above shortcomings. 
The proposed PINN approaches in this paper based on the new architecture have the following advantages:
\begin{itemize}
    \item 
    They well address the problems existing in the conventional PINN method, including: zero solution, no propagation, incorrect propagation.
    \item They use a single neural network to avoid increasing model complexity and computational cost.
    \item The overall PINN solution of the whole time domain can accurately reproduce the local solution of each time subinterval.
    \item The potential solution strictly satisfies the continuity and smoothness at the time interval nodes.
\end{itemize}

The rest of this paper is organized as follows: 
In Section 2, we first describe the conventional PINN framework for time-dependent PDEs, then give the definition of interpolation and extrapolation, and analyze the extrapolation capability of the PINN methods;
In Section 3, we construct a kind of extrapolation-driven neural network architecture by exploiting the extrapolation property;
In Section 4, we use several benchmarks to test the performance of the proposed network architecture; 
Finally, we  \textcolor{black}{provide our conclusions} in Section 5.

\section{Extrapolation property of physics-informed learning}\label{sec:Extrapolation of conventional PINNs}

In this section, we first briefly describe the conventional PINN method and then focus on discussing its extrapolation capability, which is used to improve its performance for time-dependent PDEs.

\subsection{Conventional PINN framework}\label{sec:Conventional PINN framework}

The conventional PINN method based on a fully connected neural network architecture is proposed in Ref.~\cite{RAISSI1}, which is considered as a promising framework for solving \textcolor{black}{PDEs}, and is effective for solving the evolution equation~\cite{kim2021dpm}.

Suppose time-dependent PDEs  \textcolor{black}{have} the following \textcolor{black}{general} form:
\begin{equation}\label{time-dependent PDE}
\begin{cases}
\displaystyle\frac{\partial u}{\partial t} + \mathcal{N}(u)=0,&x\in \Omega \subset \mathbb{R}^{d},t\in(0,T],\\
u(x,0)=\mathcal{I}(x),&x\in \Omega \subset \mathbb{R}^{d},\\
u(x,t)=\mathcal{B}(x, t),&x\in \partial \Omega, t\in(0,T],\\
\end{cases}
\end{equation}
where $u=u(x,t)$ is the solution on the spatio-temporal domain $\Omega \times [0, T]$, $\mathcal{N}$ denotes the differential operator with respect to the spatial variable $x$, $\mathcal{I}(x)$ and $\mathcal{B}(x,t)$ represent the initial condition (I.C.) and the Dirichlet boundary condition (B.C.), respectively. 
Boundary conditions can also be of the Neumann, Robin, periodic, and hybrid types. 

The conventional PINN method provides the prediction function $u_\theta(x,t)$, which we call the PINN solution in this paper, as an approximation of the true solution of Eq.~\eqref{time-dependent PDE}, where $\theta$ is a set of network parameters \textcolor{black}{with input $x$ and $t$} and is determined by minimizing the following loss function
\begin{equation}\label{loss_std}
\begin{cases}
\mathcal{L}(\theta;\Sigma)=w_s\mathcal{L}_s(\theta;\tau_s)+w_r\mathcal{L}_r(\theta;\tau_r),\\
\mathcal{L}_s(\theta;\tau_{s})=\frac{1}{N_0} \sum_{i=1}^{N_0}\left| u_\theta \left(x_i,0\right)-\mathcal{I}(x_i)  \right|^2 \\
\quad\quad\quad\quad + \frac{1}{N_b} \sum_{i=1}^{N_b}\left|  u_\theta \left(x_i,t_i\right)    - \mathcal{B}(x_i,t_i)  \right|^2 , \\
\mathcal{L}_r(\theta;\tau_{r})=\frac{1}{N_r} \sum_{i=1}^{N_r}\left| \frac{\displaystyle{\partial u_{\theta}}}{\displaystyle{\partial t}}(x_i,t_i) +  \mathcal{N}(u_\theta(x_i,t_i)) \right|^2.
\end{cases}
\end{equation}
In Eq.~\eqref{loss_std}, $w_s$ and $w_r$ are the weights for the two parts of the loss function,
$\Sigma$ 
denotes \textcolor{black}{a summation over} the set of training points and consists of the supervised point set $\tau_s$ and the residual point set $\textcolor{black}{\tau_r}$,
$N_0$ is the number of sample points for the I.C., $N_b$ is the number of sample points on the boundary $\partial\Omega$, and $N_r$ is the number of residual points sampled in the spatio-temporal domain.

\subsection{Extrapolation capacity analysis of PINN methods}\label{sec:Analysis of the extrapolation of PINNs}
\begin{figure}[h]
    \centering        \includegraphics[width=0.6\textwidth]{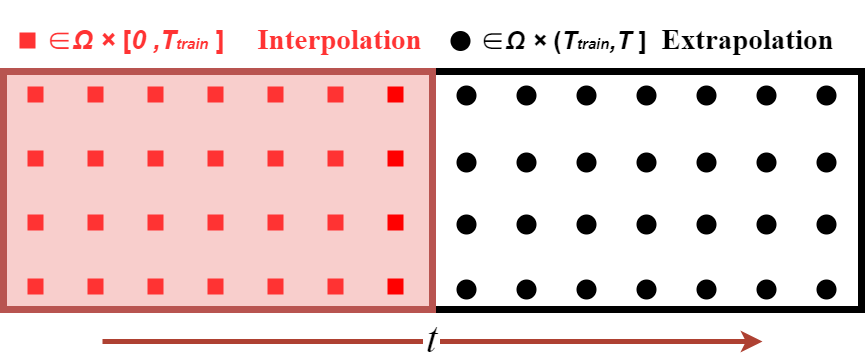}
    \caption{\textcolor{black}{Interpolation and extrapolation points.}}
    \label{fig:Interpolation and Extrapolation of PINNs}
\end{figure}

Note that the \textit{\textbf{interpolation}} and \textit{\textbf{extrapolation}} discussed in this paper are based on the time dimension.
Suppose the training domain is $\Omega \times[0, T_{train}]$.  After sufficiently sampling training points according to a certain sampling method on the training domain and optimizing the loss function~\eqref{loss_std}, we obtain the PINN solution $u_{\theta}$.
To evaluate the performance of PINNs, we usually select some test points in the training domain $\Omega \times[0, T_{train}]$. We call these test points 
\textbf{\textit{interpolation points}} here.
Previous \textcolor{black}{work} on PINN methods mainly \textcolor{black}{uses} interpolation points to verify the effectiveness of their \textcolor{black}{method}~\cite{HAGHIGHAT2021113741}.
To examine the generalization capability of the PINN solution, we need to use some test points sampled beyond the training domain, i.e., $\{(x_i,t_i)\}\in \Omega \times (T_{train}, T]$, $T>T_{train}$, we call these test points \textbf{\textit{extrapolation points}}. 
Figure~\ref{fig:Interpolation and Extrapolation of PINNs} is a schematic diagram of the interpolation and extrapolation points.

\textcolor{black}{To investigate the extrapolation capability of the PINN method}, we first consider the following
Allen-Cahn equation, which is widely used to study phase separation phenomena~\cite{C7FD00037E,ALLEN19791085},
\begin{equation}\label{llen-Cahn equation with strongly num:Allen-cahn}
\begin{cases}
u_t-0.0001u_{xx}+5u^3-5u=0, &(x,t) \in (-1,1) \times  (0,1],\\
u(x,0)=x^2 \cos{\pi x},& x\in  (-1,1),   \\
u(-1,t)=u(1,t), u_x(-1,t)=u_x(1,t), & t\in  (0,1].\\
\end{cases}
\end{equation}

\begin{figure}[ht]
    \centering    
\includegraphics[width=0.8\textwidth]{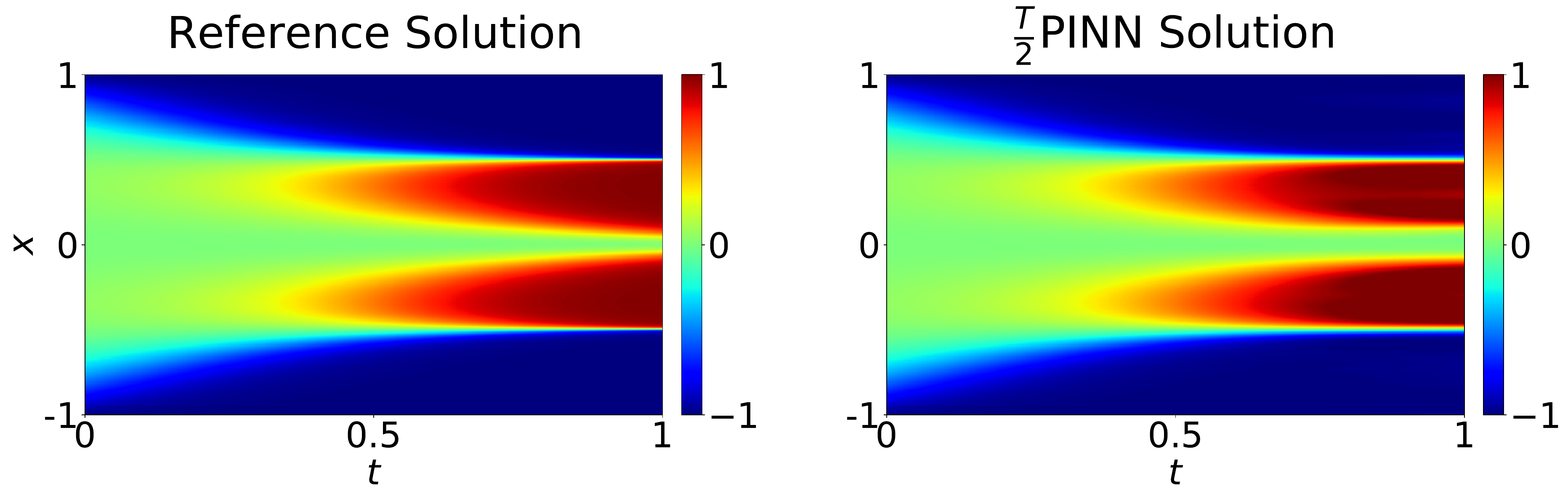}
    \caption{The reference solution of Eq.~\eqref{llen-Cahn equation with strongly num:Allen-cahn} and the \textcolor{black}{$\frac{T}{2}$PINN} solution obtained by training on the time domain $[0,0.5]$.}
    \label{fig:Interpolation and Extrapolation of PINNs in AC}
\end{figure}

\begin{figure}[h]
    \centering    
        \includegraphics[width=1\textwidth]{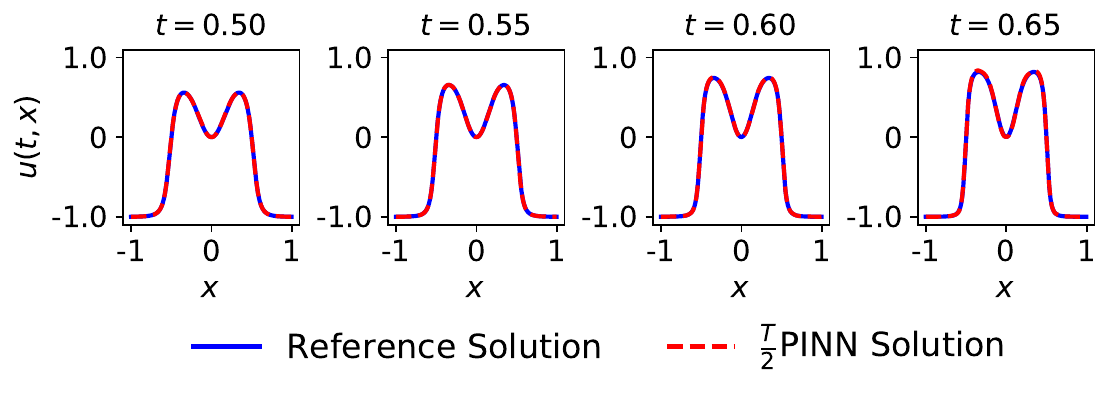}
    \caption{Comparison of the reference solution and the \textcolor{black}{$\frac{T}{2}$PINN} solution of Eq.~\eqref{llen-Cahn equation with strongly num:Allen-cahn}.}
    \label{fig:comparison of PINNs in kdv}
\end{figure}

\begin{figure}[h]
    \centering    
        \includegraphics[width=0.7\textwidth]{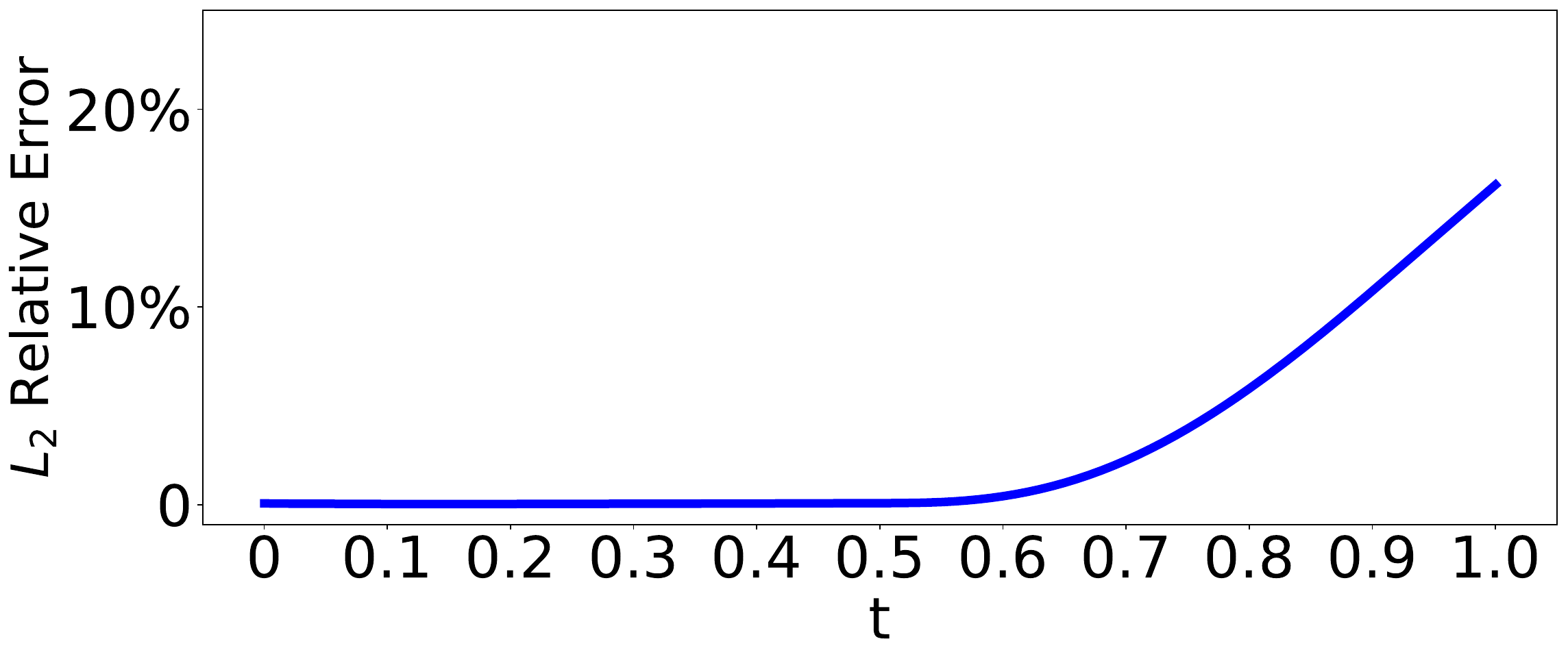}
    \caption{$L_2$ errors of $\frac{T}{2}$PINN  for Eq.~\eqref{llen-Cahn equation with strongly num:Allen-cahn}.} 
    \label{fig:comsssv}
\end{figure}

\textcolor{black}{
We train the PINNs for this PDE over half the time domain, $[-1,1] \times [0, 0.5]$, with all hyperparameters chosen as described in Section \ref{sec:Numerical Examples}. 
For simplicity and clarity, we refer to the half-domain trained PINN solution as $\frac{T}{2}$PINN.}
\textcolor{black}{Figures}~\ref{fig:Interpolation and Extrapolation of PINNs in AC} and \ref{fig:comparison of PINNs in kdv} show the contrast between the reference solution and the \textcolor{black}{$\frac{T}{2}$PINN} solution.
\textcolor{black}{We see that the $\frac{T}{2}$PINN solution also has good accuracy in predicting the values of the extrapolation points in the time domain $[0.5,1]$.
Figure~\ref{fig:comsssv} demonstrates that,   the closer $t$ to 0.5, the higher prediction accuracy.}
For this \textcolor{black}{PDE}, the PINN method has a strong extrapolation capability  \textcolor{black}{since the solution is generally smooth with respect to $t$}.

\begin{figure}[ht]
    \centering    
        \includegraphics[width=0.8\textwidth]{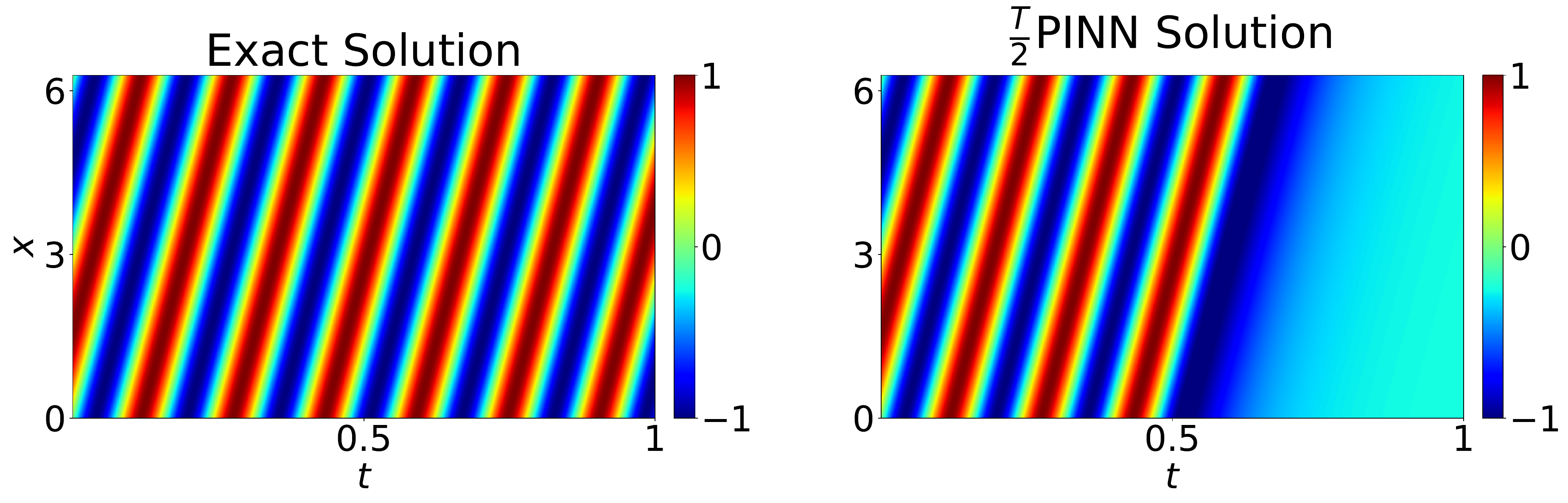}
    \caption{The reference solution of Eq.~\eqref{consider Convection and convection equation} with $\beta=40$ and the \textcolor{black}{$\frac{T}{2}$PINN} solution obtained by training on the time domain $[0,0.5]$.}
    \label{fig:Extrapolation of PINNs in Convection equation}
\end{figure}

\begin{figure}[!h!t]
    \centering    
        \includegraphics[width=0.9\textwidth]{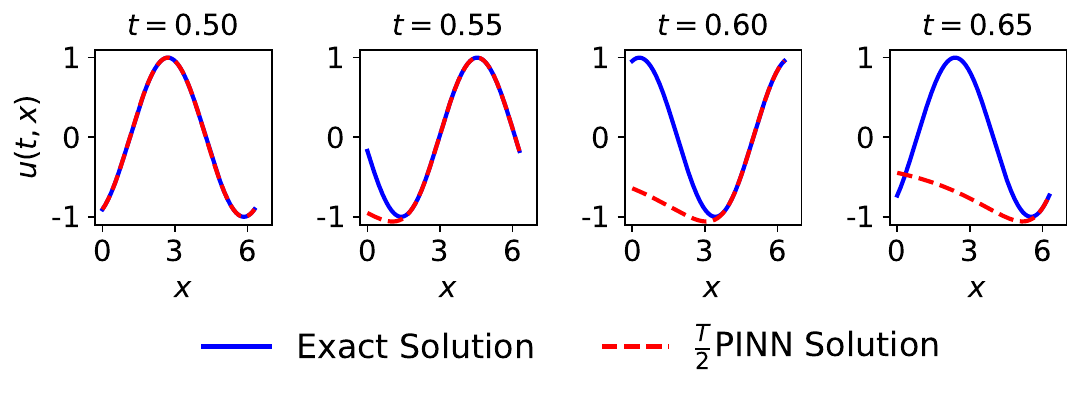}
    \caption{Comparison of the exact solution and the \textcolor{black}{$\frac{T}{2}$PINN} solution of  Eq.~\eqref{consider Convection and convection equation} with $\beta=40$.}
    \label{fig:Extrapolatin Convection equation}
\end{figure}

\begin{figure}[!h]
    \centering    
        \includegraphics[width=0.7\textwidth]{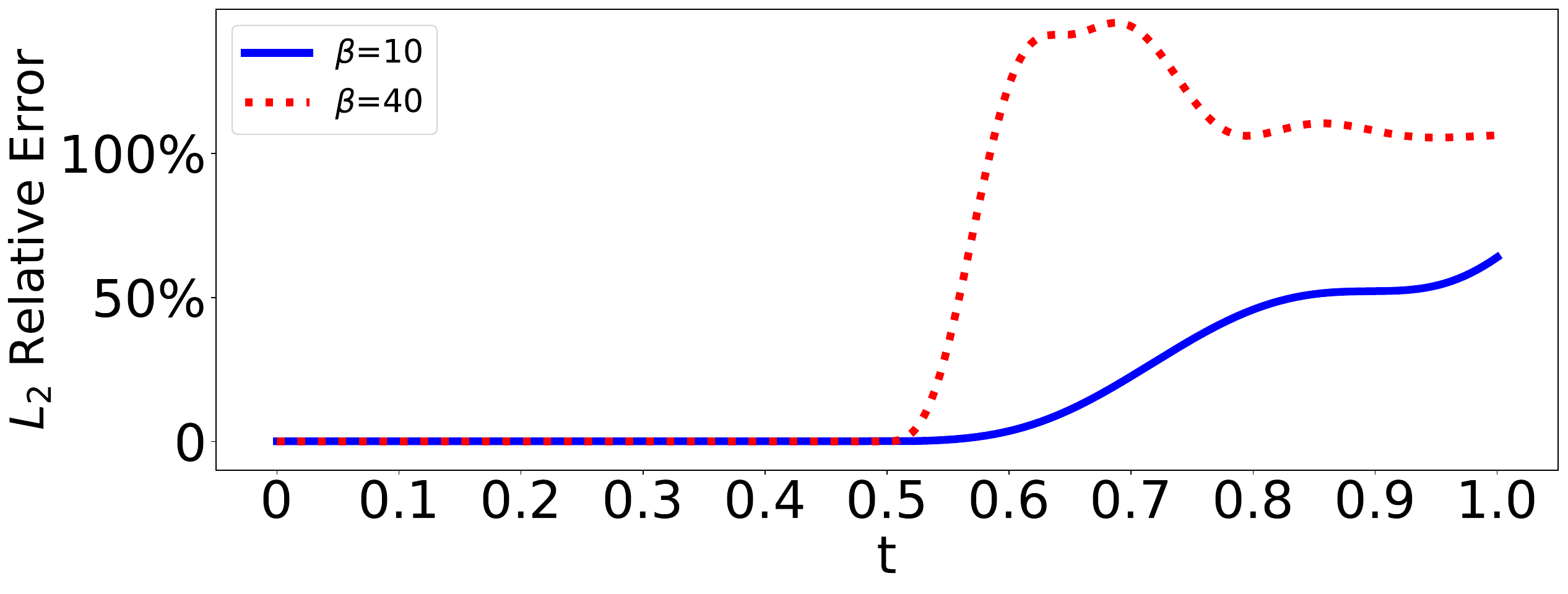}
    \caption{$L_2$ errors of $\frac{T}{2}$PINN  for Eq.~\eqref{consider Convection and convection equation} with $\beta=10$ and $\beta=40$.} 
    \label{fig:comduiliu2error}
\end{figure}

Next, we consider the
following convection equation
\begin{equation}\label{consider Convection and convection equation}
\begin{cases}
u_t+\beta u_x=0, &(x,t) \in (0,2\pi) \times  (0,1],\\
u(x,0)= \sin{ x},& x\in   (0,2\pi),   \\
u(0,t)=u(2\pi,t),& t\in  (0,1].
\end{cases}
\end{equation}
The exact solution is given by
\begin{equation}\label{Exact_Ex1}
u(x,t)=\sin{(x-\beta t)},
\end{equation}
where $\beta$ is a constant.
We let $\beta=40$, \textcolor{black}{and} train the \textcolor{black}{$\frac{T}{2}$PINN} method for this \textcolor{black}{PDE} in the computational domain $[0, 2\pi] \times[0,0.5]$.
\textcolor{black}{Figures}~\ref{fig:Extrapolation of PINNs in Convection equation}, \ref{fig:Extrapolatin Convection equation} \textcolor{black}{and~\ref{fig:comduiliu2error}}
show the contrast between the reference solution and the PINN solution.
We see that the \textcolor{black}{$\frac{T}{2}$PINN} solution trained in the time domain [0,0.5] can only make reasonable predictions for the extrapolation points very close to 0.5, but cannot make accurate predictions for the extrapolation points slightly away from 0.5. For this \textcolor{black}{PDE}, the PINN method has a weak extrapolation capability.

\begin{figure}[ht]
    \centering    
        \includegraphics[width=0.8\textwidth]{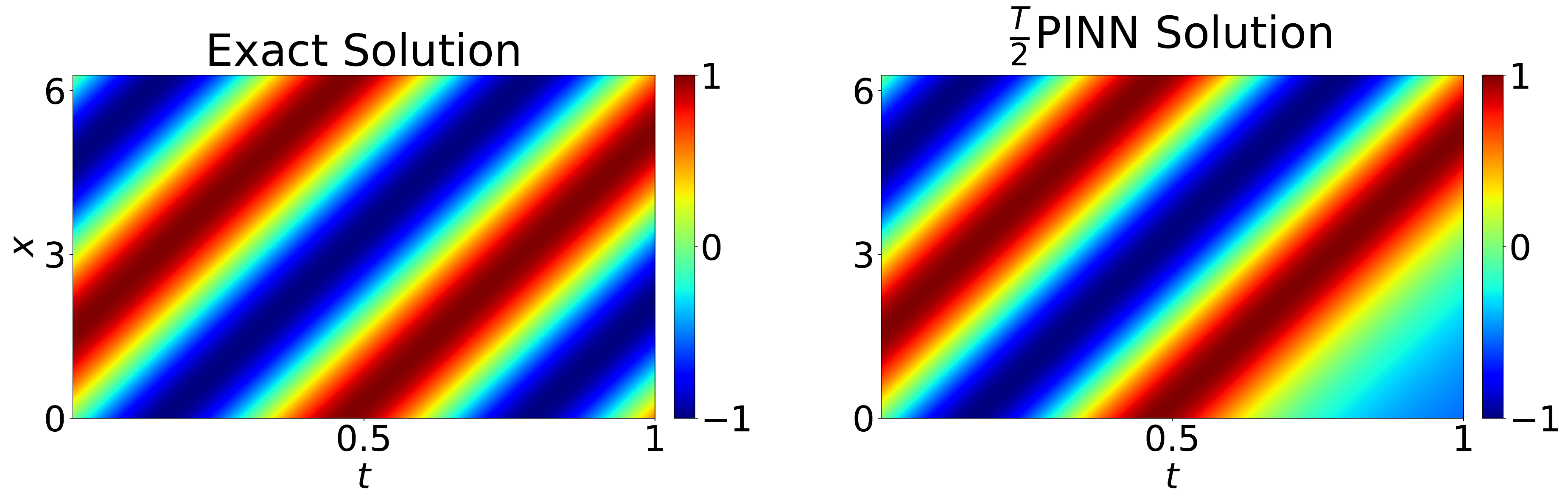}
    \caption{The reference solution of Eq.~\eqref{consider Convection and convection equation} with $\beta=10$ and the \textcolor{black}{$\frac{T}{2}$PINN} solution obtained by training on the time domain $[0,0.5]$.}
    \label{fig:10Extrapolation of PINNs in Convection equation}
\end{figure}

\begin{figure}[ht]
    \centering    
        \includegraphics[width=1.0\textwidth]{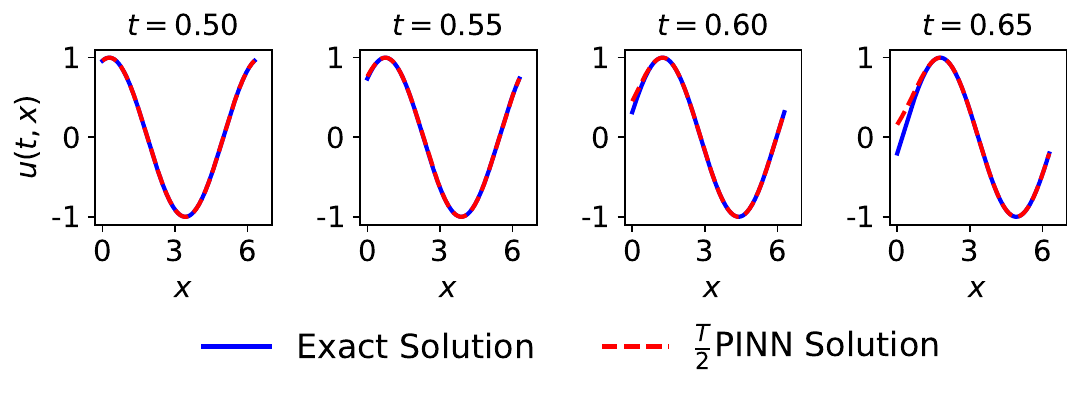}
    \caption{Comparison of the exact solution and the \textcolor{black}{$\frac{T}{2}$PINN} solution of  Eq.~\eqref{consider Convection and convection equation} with $\beta=10$.}
    \label{fig:10Extrapolatin Convection equation}
\end{figure}

We adjust the parameters in the governing equation~\eqref{consider Convection and convection equation},  
let $\beta=10$, for this case the comparison results between the exact solution and the \textcolor{black}{$\frac{T}{2}$PINN} solution are shown in \textcolor{black}{Figures}\textcolor{black}{~\ref{fig:comduiliu2error}},  ~\ref{fig:10Extrapolation of PINNs in Convection equation} 
and \ref{fig:10Extrapolatin Convection equation}.
We see that the extrapolation capability is stronger for $\beta=10$ than for $\beta=40$.
According to the exact solution, we can \textcolor{black}{see} that the value of $|u_t(x,y)|_{\beta=40}$ is greater than the value of $|u_t(x,y)|_{\beta=10}$.
A larger $|u_t|$ means that $u$ changes faster with time, and therefore the extrapolation capability is weaker.

Based on these examples, we summarize \textcolor{black}{our observations} as follows:
\begin{enumerate}
    \item \textcolor{black}{By using the governing equation as soft penalty constraints}, the PINN method has either strong or weak extrapolation capabilities~\cite{Karniadakis,HAGHIGHAT2021113741};
    \item  The prediction accuracy of the extrapolation points gradually decreases as they move away from the training domain~\cite{fesser2023understanding,bonfanti2023generalization};
    \item The extrapolation capability of the PINN method is closely related to  the \textcolor{black}{complexity of the solution}~\cite{kim2021dpm,fesser2023understanding}; 
    \item The extrapolation performance of the PINN method is weaker for the model with larger $|u_t|$ outside the training domain.
\end{enumerate}

Through analysis, we believe that for the time-dependent PDEs, there are two main reasons that make the PINN methods have some extrapolation capability:
\begin{enumerate}

\item 
The evolution equation includes the partial derivative term of the solution function $u$ with respect to time variable $t$, so $u$ is at least continuous or even smooth in the time direction; the PINN solution, as an approximation of $u$, is also continuous, which ensures that it has some accuracy for extrapolation points close to the training domain.\\
\item 
In the process of neural network training, the information of the PDE is fully exploited, so that the PINN solution approximates the true solution of the PDE in the training domain. Considering that the PDE itself holds in the entire computational domain, the PINN solution also has a certain degree of accuracy for extrapolation points in the case of slowly varying \textcolor{black}{solutions} with respect to the time variable $t$.
\end{enumerate}

\textcolor{black}{\begin{remark}
The most critical factor affecting the extrapolation capability of PINNs is the nature of the governing equation itself, which includes the type of PDE, the smoothness of the coefficients and source term in the PDE, and the well-posedness of the definite solution conditions.
The number of sampling points and the size of the neural networks do not significantly affect the extrapolation capability of PINNs (see \ref{appendixaaaaaaaa}).
\end{remark}}

\section{An extrapolation-driven neural network architecture}\label{sec:extrapolation-driven neural network architecture}
In the previous section, we discussed the extrapolation property of the PINN method. In this section, we use this property to construct new neural networks called \textbf{\textit{Extrapolation-driven Deep Neural Networks}}, abbreviated as \textbf{E-DNN}, which are helpful for solving evolution equations efficiently and accurately. 
We first give the motivation for constructing extrapolation-driven neural networks, and then present several approaches of constructing extrapolation-driven neural networks according to different characterises of PDEs.

\subsection{Motivation for constructing extrapolation-driven networks}\label{sec:motivation}

As mentioned in the introduction, \textcolor{black}{the conventional PINN faces some difficulties in solving the evolution equation along a large time dimension.}
For time-dependent PDEs, when optimizing the loss function of the conventional PINN method, since the common optimization algorithm does not take into account respecting the principle of the causality, but aims to make the value of the loss function reach the minimum, the optimization process  \textcolor{black}{is not done} in chronological order, and the positions far away from the initial time may \textcolor{black}{get} optimized earlier, which leads to its PINN solution satisfying the equation, but not matching the true solution because it is only constrained by the PDE and not affected by the initial values.  
\textcolor{black}{This resulting solution is often referred to as a trivial solution.}
As we know, the solution \textcolor{black}{that satisfies} only the PDE is usually not unique, e.g. for the heat equation $u_t-u_{xx}=f(x,t)$, if $\Tilde{u}$ satisfies the equation, then $\Tilde{u}+ax+b$ also satisfies the equation, the uniqueness is determined by combining the initial \textcolor{black}{and} boundary conditions.

To obtain high-precision solutions for the time-dependent PDEs \textcolor{black}{along a} large time \textcolor{black}{dimension}, \textcolor{black}{we should} adopt a training scheme that respects the principle of the causality.
The causality means that the solution of time-dependent PDEs at the later time is determined by the solution at the earlier time. In the training process of PINNs, \textbf{\textit{respecting the causality}} means that \textcolor{black}{the training of the later time subinterval should be based on the training results of the previous time subinterval}.
In Ref.~\cite{WANG2024116813}, the authors decompose the whole time domain into several subintervals and design a specific weight function for the residual loss term of each sub-time domain, thus achieving the goal of respecting the causal law in the training process and obtaining a high-precision solution of the evolution equation. The computational efficiency of this approach needs improvement. 
In Ref.~\cite{PENWARDEN2023112464} and Ref.~\cite{GUO2023112258}, the authors present high-precision PINN methods that respect the causal law by dividing the whole time domain into several subintervals and training them one by one. 
However, in Ref.~\cite{PENWARDEN2023112464}, the PINN solution of the whole time domain is concatenated by \textcolor{black}{all subintervals}, multiple neural networks are used, \textcolor{black}{which weakens the continuity of the solution}; 
in Ref.~\cite{GUO2023112258}, although the PINN solution is expressed by a single neural network on the whole time domain, it cannot preserve the solutions of the subintervals obtained by pre-training. 

\textcolor{black}{Motivated by the extrapolation capability of the PINN method,} we attempt to construct a novel neural network architecture based on it, and then provide an improved PINN approach from the network architecture perspective, which hopefully advances the deep learning method to solve the time-dependent PDEs efficiently and accurately. 
The main idea for constructing the new network architecture is as follows.

Suppose the whole time domain is 
$[0,T]$. 
To overcome the difficulty of training on the large time domain, we divide it into two subintervals $[0, T_p]$ and $[T_p, T]$, where $0< T_{p} < T$. 
We first train the PINNs on the first subinterval $[0,T_p]$ using the \textcolor{black}{conventional} approach, obtaining the network parameters denoted by $\theta_{p}$ and its PINN solution denoted by $u_{\theta_{p}}(x,t)$. 
\textcolor{black}{For the domain $[0,T]$, 
we can assume that the overall PINN solution}, denoted by $u_{\theta}(x,t)$, is similar to $u_{\theta_{p}}(x,t)$, so they have similar network parameters. We express their relationship as $\theta=\theta_{p}+\delta$, where $\delta$ is a set of parameters to be trained next, we call it the \textbf{\textit{correction term}}. 
To use a \textbf{\textit{single}} neural network to provide the PINN solution  $u_\theta$ for the whole time domain $[0,T]$, while requiring that $u_\theta(x,t) = u_{\theta_1}(x, t), t\in[0,T_p]$, 
we introduce a time function ${\mathcal{F}}(t)$, called \textbf{\textit{extrapolation control function}},  into the network parameters, and set 

\begin{equation}\label{ensure the T1 solution}
\theta = \theta_1 + {\mathcal{F}}(t)\cdot\delta, 
\end{equation}
thereby forming a new network architecture containing the time variable $t$. 
To ensure the smoothness of $u_\theta$ at $T_p$, we require that
\begin{equation}\label{ensure the smoothness}
\begin{cases}
\mathcal{F}(t)=0,& t\in [0, T_p],\\
\mathcal{F}(t)=1,& t \in [T,+\infty),\\
\mathcal{F}(t) \in C^1, & t\in (0,+\infty).\\
\end{cases}
\end{equation}
How to choose appropriate ${\mathcal{F}}(t)$ for equations with different properties is discussed in detail in the following sections.
By training the parameter $\delta$ on the second subinterval $[T_p, T]$, we obtain the PINN solution for the whole time domain.

\subsection{ Extrapolation-driven architecture for strong extrapolation models}\label{sec:An strong}

From the example of the Korteweg-de Vries equation in Sec.~\ref{sec:Analysis of the extrapolation of PINNs}, it can be seen that for this type of model, the PINN solution has a strong extrapolation capability, which means that the maximum norm of the correction term $||\delta||_\infty$ could be very small. Therefore, we can use any smooth function as the extrapolation control function, denoted by $\mathcal{F}_{s}(t)$.
Similar to the previous section, we suppose that the whole time domain is $[0,T]$ and divide it into two subintervals $[0, T_p]$ and $[T_p, T]$, where $0< T_{p} < T$. The control function $\mathcal{F}_{s}(t)$ should meet the following conditions:

\begin{enumerate}
  
    \item $\mathcal{F}_{s}(t)=0, t\in [0,T_{p}].$\\ 
    This condition is very important, it serves two purposes. First, it introduces the training results of the first subinterval into the second subinterval as equivalent initial conditions; second, it prevents re-optimization of the first subinterval, thus preserving the accuracy of $[0, T_p]$, since the PINN method often has higher accuracy at smaller intervals~\cite{GUO2023112258,10191822}. 

    \item  $\mathcal{F}_s(t) =1, t\in [T,+\infty).$\\
   This condition provides an upper bound for the control function $\mathcal{F}_{s}(t)$, which on the one hand helps to optimize the correction term on the subinterval $[T_p, T]$, and on the other hand can provide fixed initial values of the network parameters for the subinterval after time T, thus making our method scalable.

    \item  $\mathcal{F}_s(t) \in C^1, t\in (0, +\infty).$\\
This condition makes the network parameters change continuously from the network parameter $\theta_p$ of the subinterval $[0,T_p]$ to the network parameter $\theta$ of the whole time domain $[0,T]$, ensuring the continuity of the PINN solution and its first partial derivative with respect to $t$ at the interval node $T_p$.
If the evolution equation contains $u_{tt}$, we require that
$\mathcal{F}_{s}(t)\in C^2$.
    
    \item  $\mathcal{F}_{s}'(t) \geq 0, t\in (T_{p},T).$\\
This condition is optional, not mandatory. We use it to cater for the second phenomenon observed in Sec.~\ref{sec:Analysis of the extrapolation of PINNs} that the extrapolation capability gradually weakens as time goes on.
It is designed with the idea that the further away the time is from $T_p$, the more the network parameters should change.

\end{enumerate}

\begin{figure}[ht]
    \centering    
    \includegraphics[width=0.5\textwidth]{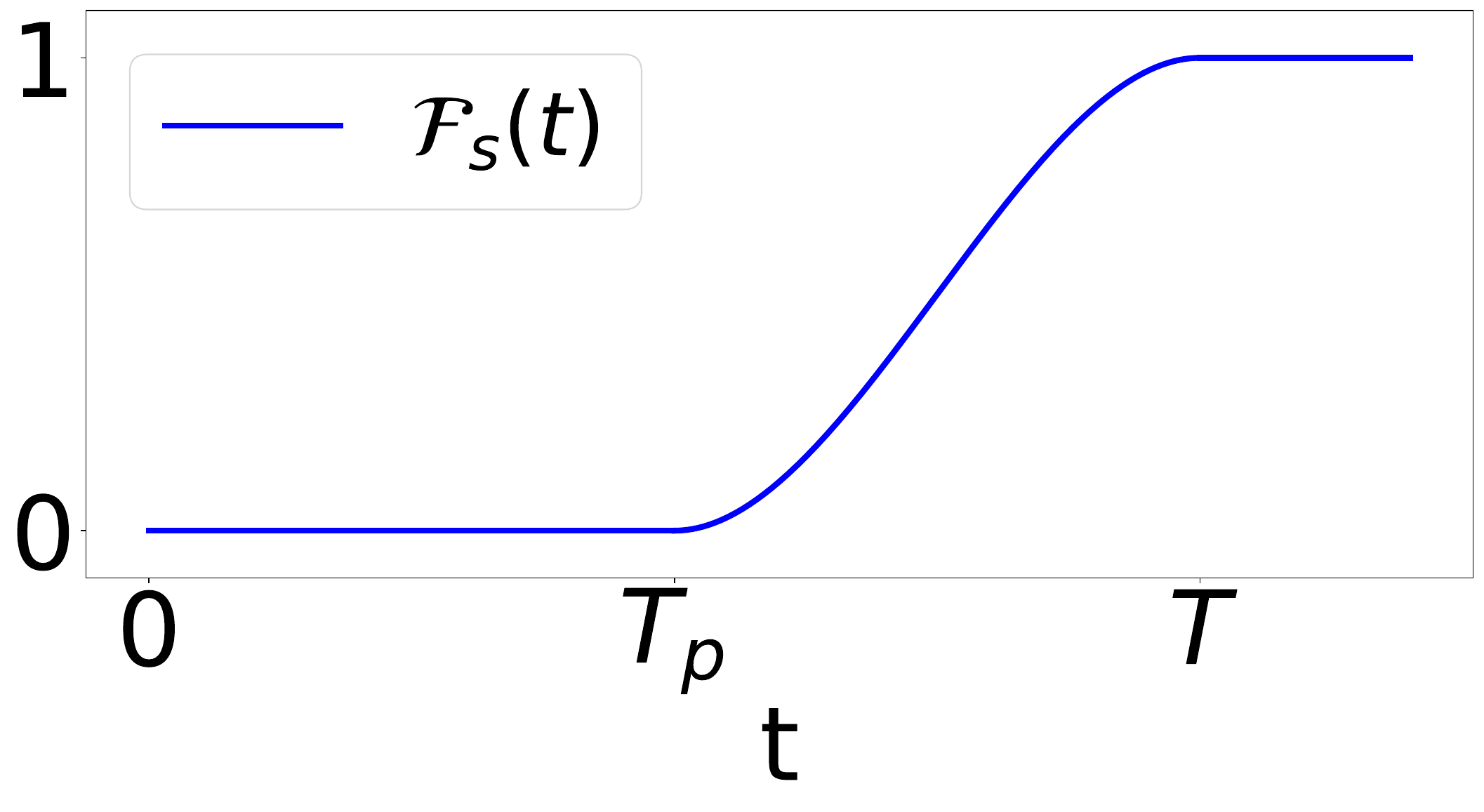}
    \caption{Graph of the extrapolation control function~\eqref{fsttt}.}
    \label{fig:Incorporating the function fst}
\end{figure}

Based on the above conditions, we provide a specific extrapolation-driven function as follows:

\begin{equation}\label{fsttt}
\mathcal{F}_{s}(t) = \begin{cases}
0, & t\in [0, T_{p}),\\
-2(\frac{t-T_p}{T-T_p})^3+3(\frac{t-T_p}{T-T_p})^2, &t\in [T_{p},T),\\
1, & t\in [T,+\infty).
\end{cases}
\end{equation}
Figure~\ref{fig:Incorporating the function fst} provides  a visual representation of $\mathcal{F}_{s}(t)$. 

\textcolor{black}{\begin{remark}
In the E-DNN architecture, the choice of the extrapolation control function is not unique.
For any function satisfying the above conditions, function \eqref{fsttt} can be viewed as an approximation of it by a third-order Hermite interpolation polynomial.
We investigate the impact of different control functions on our approach, and the numerical results are presented in the \ref{appendixBBBBBBBB}.
\end{remark}}

\begin{figure}[ht]
    \centering    
    \includegraphics[width=0.8\textwidth]{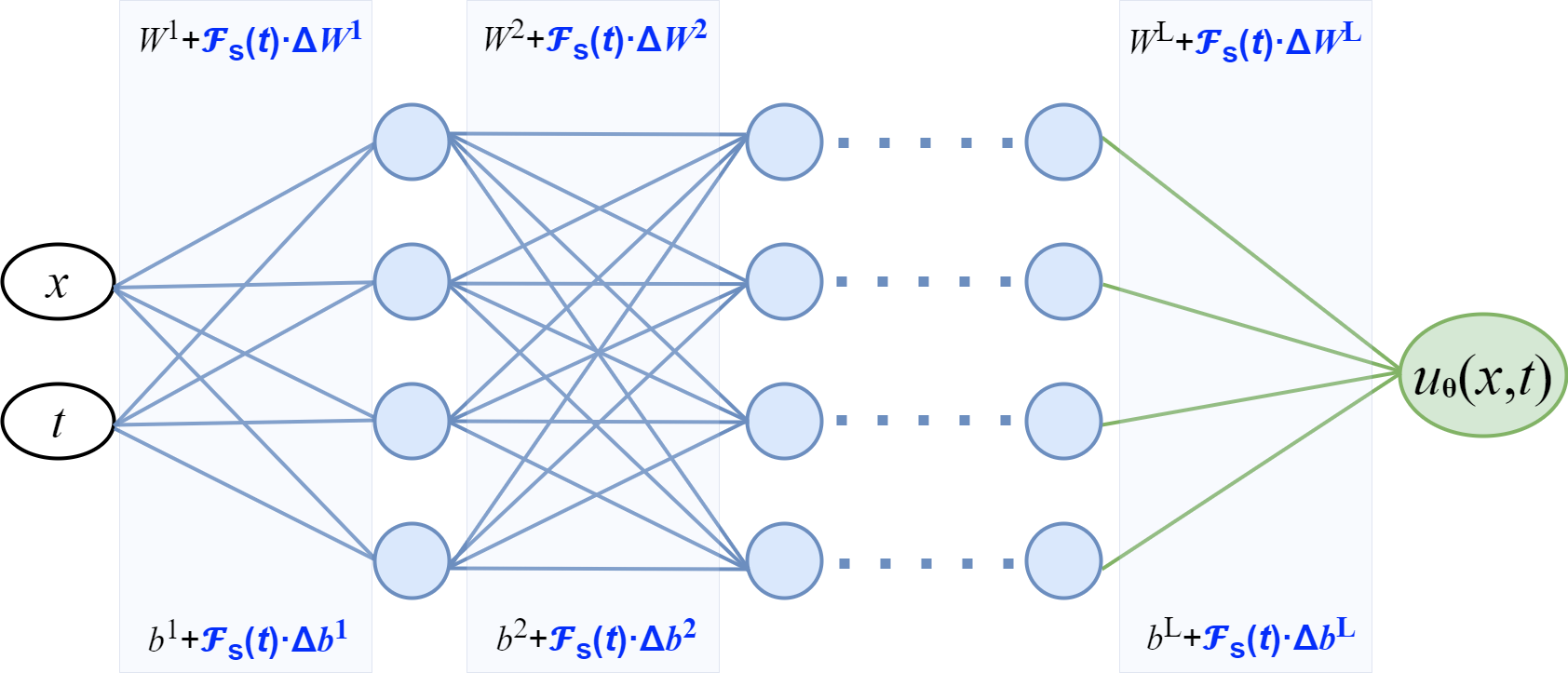}
    \caption{The E-DNN architecture.}
    \label{fig:SENN}
\end{figure}

To obtain the overall PINN solution of PDEs, our training scheme is as follows: First, we use the usual feed-forward fully connected deep neural networks on the subinterval $[0,T_{p}]$ and obtain their network parameters, including the weights $W$ and biases $b$; then, we use the extrapolation-driven deep neural network architecture, abbreviated as \textbf{\textit{E-DNN architecture}}, to obtain the PINN solution of the time domain $[0,T]$.
The E-DNN architecture for strong extrapolation models is shown in Figure~\ref{fig:SENN}.

The detailed feed-forward pass process of the E-DNN architecture is formulated as follows:

\begin{equation}\label{the detailed feed-forward pass of the S-EDNN}
\begin{cases}
Z^1 = \varnothing((W^1+ \mathcal{F}_{s}(t)\Delta W^1)   X+     b^1+\mathcal{F}_{s}(t)\Delta b^1),&\\
Z^k = \varnothing((W^k+\mathcal{F}_{s}(t)\Delta W^k)Z^{k-1}+b^k+\mathcal{F}_{s}(t)\Delta b^k),& k=2, \cdots, L-1,\\
u_{\theta}(x,t) = (W^L+\mathcal{F}_{s}(t)\Delta W^L)Z^{L-1}+b^L+\mathcal{F}_{s}(t)\Delta b^{L},\\
\end{cases}
\end{equation}
where $L$ is the number of  hidden layers of the neural networks, $X=(x,t)$ denotes the vector of input variables, $\varnothing$ is a nonlinear activation function,
$\{(W^k,b^k)_{k=1}^{L}\}$ is a set of neural network parameters obtained by prior training on the subinterval $[0, T_p]$, which are fixed in this E-DNN architecture, and $\{ (\Delta W^k, \Delta b^k)_{k=1}^{L} \}$ are the parameters of the E-DNN architecture to be trained.
\textcolor{black}{Note that E-DNN  does not increase the size of the network.}
\textcolor{black}{\begin{remark}
The E-DNN architecture exploits the extrapolation capability of PINNs on the domain $\Omega \times [T_{train}, T]$, which is central to the design of the extrapolation control function for sequential learning under the time-domain decomposition framework. 
Note that $\mathcal{F}(t)$ and $\delta$ are new parameters of the neural network in the interval $[T_{p}, T] \times \Omega$ that need to be determined after training. In this sense, the E-DNN approach aims to improve the interpolation performance of PINNs over $[0,T]$ rather than to predict solutions for $t>T$.
\end{remark}}

\subsection{Extrapolation-driven architecture for weak extrapolation models}\label{sec:An weak Extrapolation-Driven Neural Network Architecture extrapolation-driven neural network architecture}

For the case of the convection equation with $\beta=40$ in Sec.~\ref{sec:Analysis of the extrapolation of PINNs}, its PINN solution has a weak extrapolation capability, which means that the correction term should be large and we have to use a \textbf{rapidly increasing function} as the extrapolation control function, otherwise we cannot obtain a highly accurate overall PINN solution.
We denote the extrapolation control function for this type of model by $\mathcal{F}_{w}(t)$, it should also satisfy the conditions mentioned in Sec.~\ref{sec:An strong}.

\begin{figure}[!h]
    \centering    \includegraphics[width=0.5\textwidth]{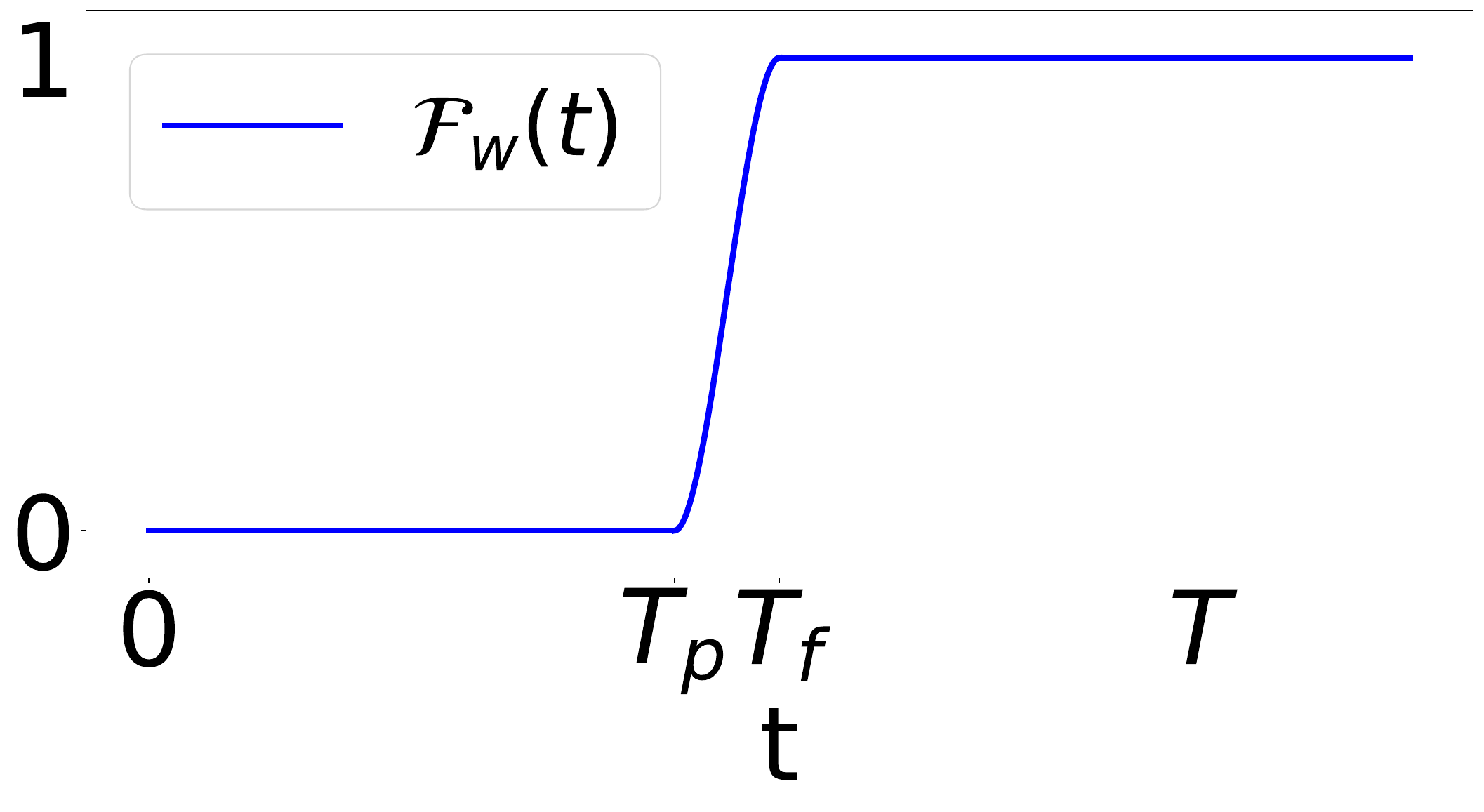}
    \caption{Graph of the extrapolation control function~\eqref{fwwwwt} with $\mathcal{M}=5$.}
    \label{fig:Incorporating the function fwt}
\end{figure}

By adding the parameter $T_f$, we provide an adjustable extrapolation control function as follows:

\begin{equation}\label{fwwwwt}
\mathcal{F}_{w}(t) = \begin{cases}
0, & t\in [0, T_{p}),\\
-2(\frac{t-T_p}{T_f-T_p})^3+3(\frac{t-T_p}{T_f-T_p})^2, & t\in [T_{p},T_f),\\
1, & t\in [T_{f},+\infty),
\end{cases}
\end{equation}
where $T_f=T_p+\frac{T-T_p}{\mathcal{M}}$, ${\mathcal{M}}$ is a positive integer constant.
Figure~\ref{fig:Incorporating the function fwt} shows the appearance of this function with $\mathcal{M}=5$.
We can see that it increases rapidly from 0 to 1, which is suitable for models with weak extrapolation capability.
For models with a weaker extrapolation capability, the extrapolatable zone becomes narrower, we need to set $\mathcal{M}$ to a larger value. 

The E-DNN architecture for the weak extrapolation model is similar to that shown in Figure~\ref{fig:SENN}, and its feed-forward pass process is similar to Eq.~\eqref{the detailed feed-forward pass of the S-EDNN}, only replacing the control function $\mathcal{F}_{s}(t)$ with $\mathcal{F}_{w}(t)$.

\subsection{Model-adaptive extrapolation-driven architecture}\label{sec:Model-adaptive}

When we apply the E-DNN architecture to solve PDEs, it is usually difficult to determine the extrapolation capability of the PINN solution for any model. Even if we can evaluate the extrapolation capability for some models according to the property of PDEs, we cannot easily determine the optimal value of $T_f$ in Eq.~\eqref{fwwwwt}.
To address these issues, we take advantage of the feature learning capabilities of neural networks~\cite{sdsads,726791} and take $T_f$ as a trainable parameter in the extrapolation control function.
We denote the extrapolation control function containing the trainable parameter $T_f$ by $\mathcal{F}_{a}(t)$, which is given by the following formula.

\begin{equation}\label{faaaat}
\mathcal{F}_{a}(t) = \begin{cases}
0, & t\in [0, T_{p}),\\
-2(\frac{t-T_p}{T_f-T_p})^3+3(\frac{t-T_p}{T_f-T_p})^2, & t\in [T_{p},T_f),\\
1, & t\in [T_{f},+\infty),
\end{cases}
\end{equation}
where $T_f \in (T_p,T]$, it is initialized with $T_p+\frac{T-T_p}{2}$ in our numerical experiments.

Similarly, we require that $\mathcal{F}_{a}(t)$ satisfies the conditions mentioned in Sec.~\ref{sec:An strong}.
By replacing the control function $\mathcal{F}_{s}(t)$ with $\mathcal{F}_{a}(t)$ in Figure~\ref{fig:SENN}, we get the model-adaptive extrapolation-driven architecture.
Note that in addition to $\{(\Delta W^k, \Delta b^k)_{k=1}^{L}\}$ being the parameters to be trained, $T_f$ is also a key parameter to be trained for this architecture.

\subsection{Two comments on PINN approaches using the E-DNN architecture }\label{sec:Two remarks}

\begin{enumerate}
\item \textbf{Training is completely independent between subintervals.}\\
As mentioned above, when using the E-DNN architecture to extend the PINN solution from the subinterval $[0, T_p]$ to $[0, T]$, the parameters we need to train are $\{(\Delta W^k, \Delta b^k)_{k=1}^{L}\}$ or $\{(\Delta W^k, \Delta b^k)_{k=1}^{L},T_f\} $. Training these parameters requires only sample points in the subinterval $(T_p, T]$ and not sample points in the subinterval $[0, T_p]$, so the training process is completely independent of the subinterval $[0, T_p]$, which makes our approach computationally efficient.
\item 
\textbf{Our approach is inherently scalable.}\\
After obtaining the PINN solution of the time interval $[0, T]$ using the above E-DNN architecture, if we further want to obtain the PINN solution of the extended subinterval $[T, \widetilde{T}], \widetilde{T} > T$, we just need to treat the interval $[0, T]$ as $[0, T_P]$ and the subinterval $[T, \widetilde{T}]$ as $[T_p, T]$, then we can get it using the same methods discussed in the above subsections. The E-DNN architecture of the subinterval $[T, \widetilde{T}]$ is given below.

We use $\{\widetilde{W}, \widetilde{b}\}$ to denote the network parameters of the interval $[0,T]$. 
Since they are trained by the E-DNN, they have the following expressions: 
\begin{equation}\label{tdwtdb }
\begin{cases}
\widetilde{W} = W + \mathcal{F}(t) \Delta W,\\
\widetilde{b} = b + \mathcal{F}(t) \Delta b.\\
\end{cases}
\end{equation}
According to the conditions mentioned in Sec.~\ref{sec:An strong}, $\mathcal{F}(t)=1, t \in [T, +\infty)$.

The feed-forward pass process of the E-DNN for $[T,\widetilde{T}]$ is as follows:
\begin{equation}\label{the expanded E-DNN architecture }
\begin{cases}
Z^1 = \varnothing((\widetilde{W}^1+ \widetilde{\mathcal{F}}(t)\Delta \widetilde{W}^1)   X+     \widetilde{b}^1+\widetilde{\mathcal{F}}(t)\Delta \widetilde{b}^1),&\\
Z^k = \varnothing((\widetilde{W}^k+\widetilde{\mathcal{F}}(t)\Delta \widetilde{W}^k)Z^{k-1}+\widetilde{b}^k+\widetilde{\mathcal{F}}(t)\Delta \widetilde{b}^k),& k=2, \cdots, L-1,\\
u_{\theta}(x,t) = (\widetilde{W}^L+\widetilde{\mathcal{F}}(t)\Delta \widetilde{W}^L)Z^{L-1}+\widetilde{b}^L+\widetilde{\mathcal{F}}(t)\Delta \Tilde{b}^{L}.\\
\end{cases}
\end{equation}
$\widetilde{\mathcal{F}}(t)$ is the extrapolation control function for the subinterval $[T, \widetilde{T}]$, one of its formulas is  
\begin{equation}\label{fsssttt}
\widetilde{\mathcal{F}}_{s}(t) = \begin{cases}
0, & t\in [0, T),\\
-2(\frac{t-T}{\widetilde{T}-T})^3+3(\frac{t-T}{\widetilde{T}-T})^2, &t\in [T,\widetilde{T}),\\
1, & t\in [\widetilde{T},+\infty).
\end{cases}
\end{equation}

Note that $\{(\Delta \widetilde{W}^k, \Delta \widetilde{b}^k)_{k=1}^{L}\}$ are the parameters to be trained in the subinterval $[T,\widetilde{T}]$.
\end{enumerate}
\textcolor{black}{
\begin{remark}
The E-DNN architecture allows us to divide a large time domain into multiple subintervals and solve the time-dependent PDEs sequentially, in a chronological order. 
Generally, training on a smaller domain yields more accurate predictions. 
By respecting causality, physical information from previous subintervals  is accurately transmitted to subsequent subintervals, effectively addressing challenges associated with large time domains, such as zero solution, no propagation, and incorrect propagation.
\end{remark}
}

\section{Numerical examples}\label{sec:Numerical Examples}
In this section, we demonstrate the performance of the proposed architectures by solving several challenging time-dependent benchmarks, including the  Allen-Cahn equation, the convection equation, and the Korteweg-de Vries equation.
We use the deep learning framework TensorFlow (version 1.13.1) to implement all experiments. 
Except for \textcolor{black}{certain cases that we will mention}, the data type is \texttt{float64} and the activation function is \texttt{tanh} for all examples.
We use the Adam optimizer~\cite{kingma2017adam} to run 5000 epochs and then switch to the L-BFGS optimizer~\cite{LLLBFGS} until convergence.
All parameters and termination criteria of the L-BFGS optimizer are considered as suggested in Ref.~\cite{LLLBFGS}.
We choose $T_p=\frac{T}{2}$ to divide the whole time domain $[0,T]$ into two subintervals.
Table~\ref{table1_sssex1} lists the PINN versions considered in this section.

\begin{table}[!h]
		\setlength{\abovecaptionskip}{0cm}
		\setlength{\belowcaptionskip}{0.2cm}
\caption{PINN versions considered in the numerical examples.}
\label{table1_sssex1}
\centering
\begin{adjustbox}{max width=\textwidth}
\begin{tabular}{c|c} 
\hline 
PINN versions &Description\\
\hline
$\frac{T}{2}$PINN & PINN solution of $[0,\frac{T}{2}]$\\
\hline
sE-PINN & Using E-DNN architecture with ${\mathcal{F}}_s(t)$ \ref{fsttt} \\
\hline
wE-PINN & Using E-DNN architecture with ${\mathcal{F}}_w(t)$ \ref{fwwwwt} \\
\hline
aE-PINN & Using E-DNN architecture with ${\mathcal{F}}_a(t)$ \ref{faaaat} \\
\hline
\end{tabular}
\end{adjustbox}
\end{table}

The $L_2$ relative error norm is used to evaluate the accuracy of different PINN versions, which is defined as follows:
\begin{equation}\label{L2norm}
\left\|\epsilon\right\|_2=\frac{\sqrt{\sum_{i=1}^N\left|u_\theta(x_i,t_i)-u(x_i,t_i)\right|^2}}{\sqrt{\sum_{i=1}^N\left|u(x_i,t_i)\right|^2}},
\end{equation}
where $N$ is the number of test points, $u(x_i,t_i)$ is the exact or reference solution of point $(x_i,t_i)$, $u_\theta(x_i,t_i)$ is the PINN solution of point $(x_i,t_i)$.

Note that we choose the Xavier scheme~\cite{XuC} to initialize the network parameters for the first subinterval $[0, T_p]$,
while we employ zero-initialization for the correction term of the second subinterval $[T_p,T]$. \textcolor{black}{In addition, 
for a fair comparison, all methods within each subsection use the same network parameters, including the same number of hidden layers and neurons, and the same width of the NN.}

\subsection{Allen-Cahn equation with highly nonlinear source term}\label{sec:Allen-Cahn equation}
Let us revisit the Allen-Cahn equation~\eqref{llen-Cahn equation with strongly num:Allen-cahn}. 
In the literature~\cite{PENWARDEN2023112464,WANG2024116813,CiCP-29-930,MATTEY2022114474,GUO2023112258,MCCLENNY2023111722}, the authors all use this model as a benchmark to test the performance of their methods. 
As mentioned in Ref.~\cite{PENWARDEN2023112464}, the larger the time domain, the more difficult it is to solve this equation using the conventional PINN method.
To test the performance of our approaches, we use sE-PINN and aE-PINN to solve this equation.

We set 4 hidden layers for the E-DNN architecture, each layer with 50 neurons, the number of initial condition points $N_0=200$, the number of periodic boundary condition points $N_b=800$, and the number of residual points $N_r=10000$. 
For practical consideration, we construct a Fourier feature embedding of the input to enforce periodic boundary conditions~\cite{WANG2024116813,DONG2021110242}.

\begin{figure}[!h]
    \centering    
    \includegraphics[width=0.99\textwidth]{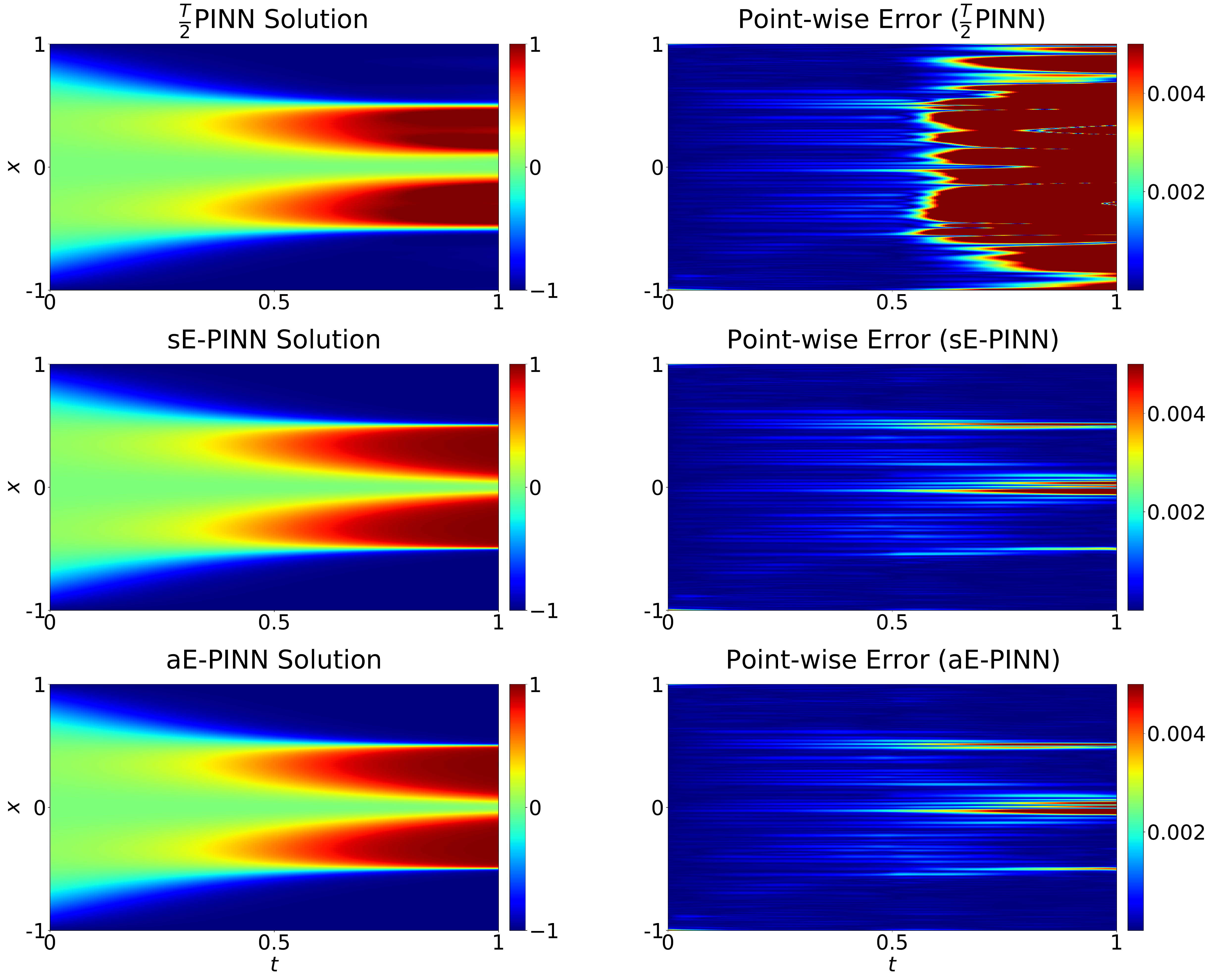}
    \caption{Numerical solutions and point-wise errors of different PINN versions for Eq.~\eqref{llen-Cahn equation with strongly num:Allen-cahn}.}
    \label{fig:ac3 equation conventional pinn and s-ednn}
\end{figure}

Figure~\ref{fig:ac3 equation conventional pinn and s-ednn} shows the images of solutions and point-wise errors for different PINN versions(float32). 
If we focus only on the subinterval $[0,0.5]$, we can see that sE-PINN and aE-PINN have the same point-wise errors as that of the $\frac{T}{2}$PINN, which means that the E-DNN architecture successfully reproduces the solution of the previous subintervals.

\begin{figure}[!h]
    \centering    
    \includegraphics[width=0.95\textwidth]{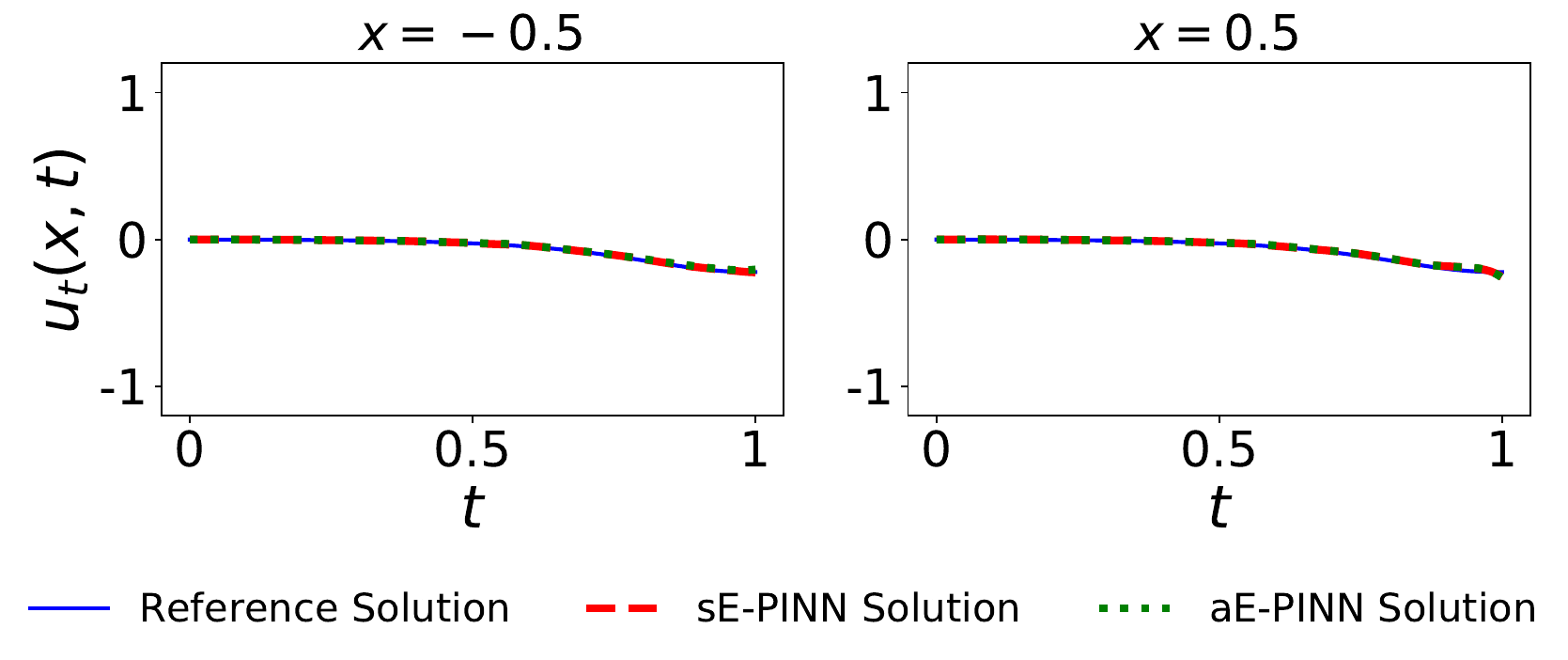}
    \caption{Comparison of $u_t(x,t)|_{x = -0.5,0.5}$  between different PINN versions for Eq.~\eqref{llen-Cahn equation with strongly num:Allen-cahn}.}
    \label{ACACACACthe continuity and smoothness at the time-split position}
\end{figure}

Since sE-PINN and aE-PINN are based on the single network, their solutions are continuous and smooth at the time interval nodes. \textcolor{black}{Figure~\ref{ACACACACthe continuity and smoothness at the time-split position} shows that the partial derivatives $u_t(x,t)$ of the PINN solutions(float32) obtained by E-DNN approximate well to  those of the reference solution at  $x=-0.5$ and  $x=0.5$.}

\begin{figure}[!ht]
    \centering    
    \includegraphics[width=0.7\textwidth]{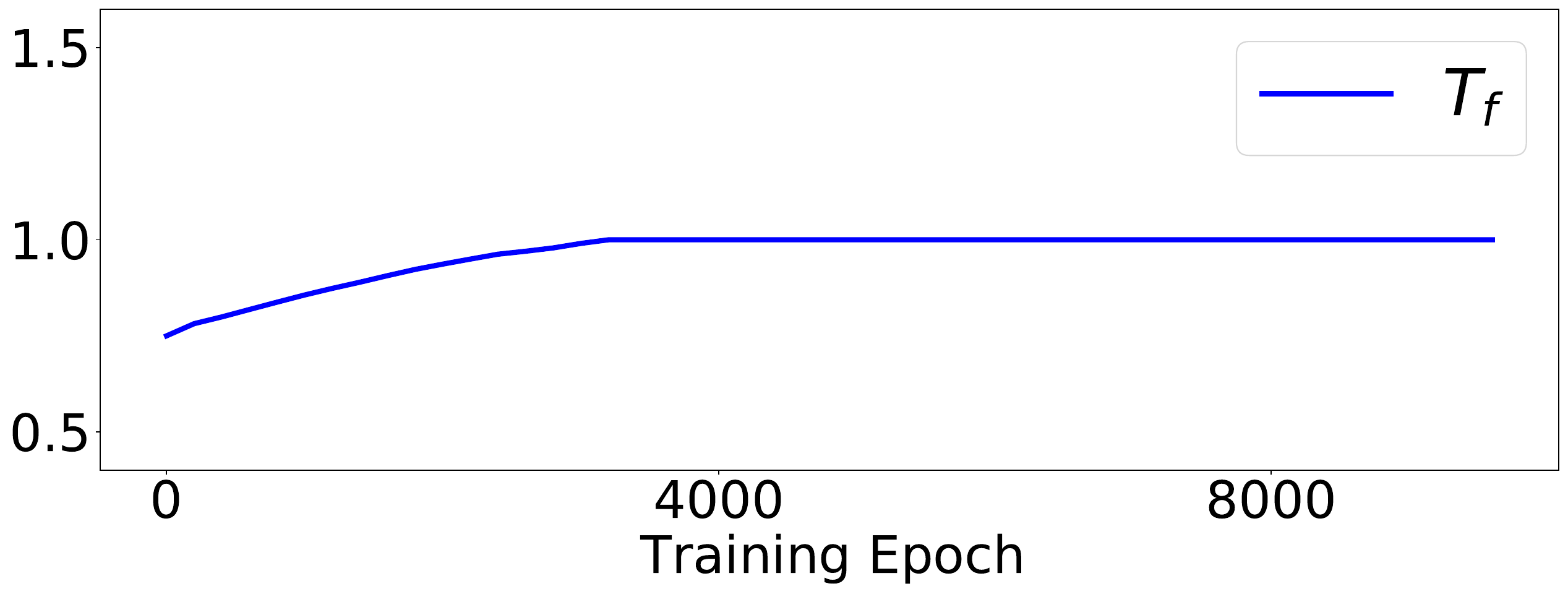}
    \caption{\textcolor{black}{Values of the trainable parameter $T_{f}$ in aE-PINN method for Eq.~\eqref{llen-Cahn equation with strongly num:Allen-cahn}}.}
    \label{The learnable parameters consider ACACACACACA equation}
\end{figure}

Figure~\ref{The learnable parameters consider ACACACACACA equation} shows the \textcolor{black}{values} of the trainable parameter $T_f$ during the optimization process when the adaptive control function $\mathcal{F}_a(t)$ is chosen. We can see that it is eventually optimized from the initial value of 0.75 to 1, which demonstrates the ability of the adaptive approach, it correctly identifies this model as one with strong extrapolation capabilities. 

\begin{figure}[!ht]
    \centering    
    \includegraphics[width=0.95\textwidth]{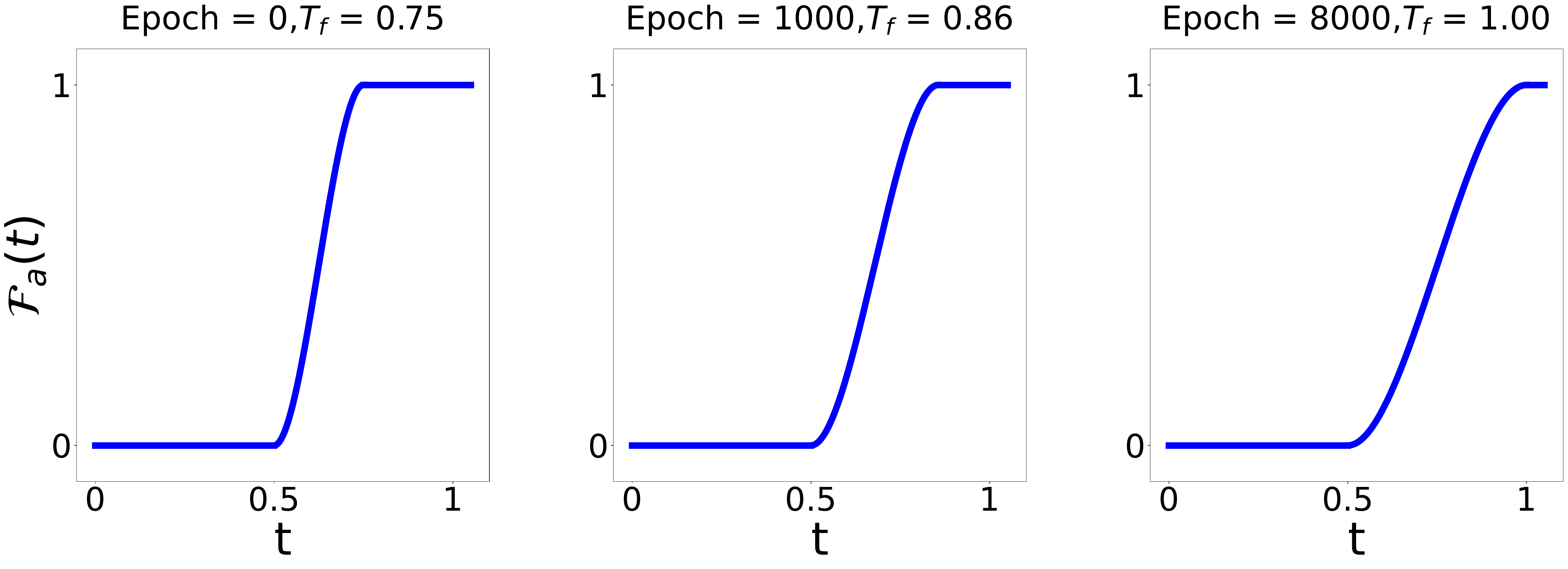}
    \caption{Optimization process of adaptive control function $\mathcal{F}_a(t)$ for Eq.~\eqref{llen-Cahn equation with strongly num:Allen-cahn}.}
    \label{2The learnable  parameters consider ACACACACACA equation}
\end{figure}

Figure~\ref{2The learnable  parameters consider ACACACACACA equation} shows the optimization process of the adaptive control function  $\mathcal{F}_a(t)$, it gradually becomes smoother from its initial steeper shape to match the strong extrapolation property of this model.

Table~\ref{table Allen-Cahn equation with strongly nonlinear source term}  lists the \textcolor{black}{accuracy (given by the $L_2$ relative error) and training time for the} various methods. Considering that \textcolor{black}{the float32 data type is often used}
in the literature, we also \textcolor{black}{use} the float32 data type for this test.
Compared to the other methods proposed in \textcolor{black}{recent years}, our methods provide more accurate predictions and show better efficiency. 
Compared to the sE-PINN method, the aE-PINN method introduces an additional trainable parameter ($T_f$), which slightly increases the training time.

\begin{table}[h]
\color{black} 
\setlength{\abovecaptionskip}{0cm}
		\setlength{\belowcaptionskip}{0.2cm}
\caption{The performance of different methods for solving Eq.~\eqref{llen-Cahn equation with strongly num:Allen-cahn}.}
\label{table Allen-Cahn equation with strongly nonlinear source term}
\centering
\begin{adjustbox}{max width=\textwidth}
\begin{tabular}{c|c|c}
\hline 
Method  &  $\left\|\epsilon\right\|_2$  & \textcolor{black}{Training Time}\\
\hline
SA-PINN~\cite{MCCLENNY2023111722} & $7.97\times 10^{-2}$ & 53.68 min\\ 
\hline
bc-PINN~\cite{MATTEY2022114474}& $2.14\times 10^{-2}$ & 105.52 min\\ 
\hline
w-s PINN~\cite{PENWARDEN2023112464}& $1.85\times 10^{-2}$ & 29.62 min\\ 
\hline
Casual Weight~\cite{WANG2024116813}& $1.71\times 10^{-2}$ & 12h++\\ 
\hline
PT-PINN~\cite{GUO2023112258}& $8.43\times 10^{-3}$ & 9.91 min\\ 
\hline
Adaptive Resampling~\cite{CiCP-29-930}& $7.39 \times 10^{-3}$ & 21.02 min\\ 
\hline
sE-PINN (float32) this work& $1.24\times 10^{-3}$ & 13.87 min\\ 
\hline
aE-PINN (float32) this work& $1.25\times 10^{-3}$ & 14.21 min\\  
\hline
aE-PINN (2 subintervals) this work& $9.36\times 10^{-4}$ & 30.25 min\\ 
\hline
aE-PINN (5 subintervals) this work& $7.86\times 10^{-4}$ & 24.37 min\\ 
\hline
\end{tabular}
\end{adjustbox}
\end{table}



\begin{figure}[!h]
    \centering    
\includegraphics[width=0.99\textwidth]{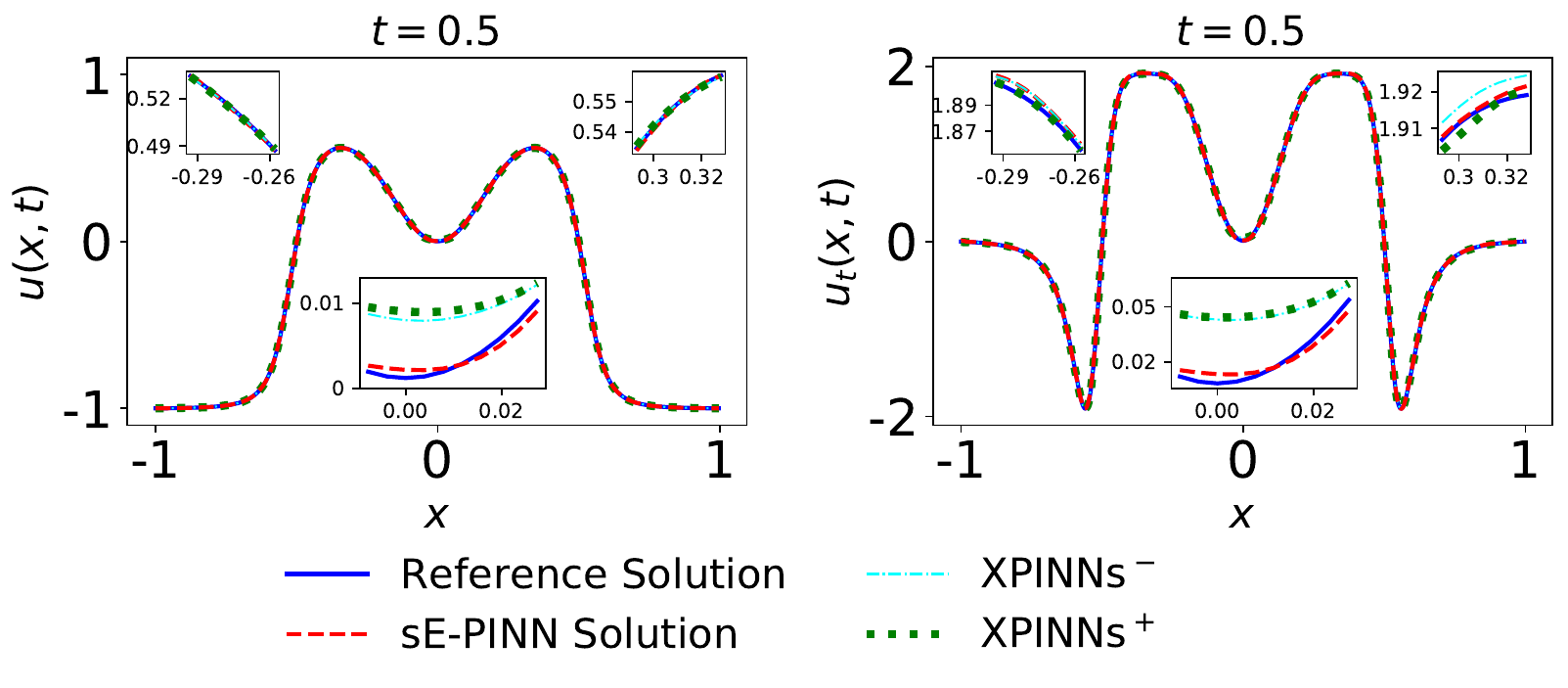}
    \caption{Comparison of $u(x,t)|_{t = 0.5}$ and  $u_t(x,t)|_{t = 0.5}$ between different PINN methods for Eq.~\eqref{llen-Cahn equation with strongly num:Allen-cahn}.} 
    \label{fig:xpinssandothers}
\end{figure}

\textcolor{black}{
We divide the computational domain into two parts, $[-1,1] \times  [0,0.5]$ and $[-1,1] \times  [0.5,1]$, and apply the XPINNs method to solve Eq.~\eqref{llen-Cahn equation with strongly num:Allen-cahn}.
Table~\ref{HEXPINNSDEBIJIAO} and  Figure~\ref{fig:xpinssandothers}  provide a detailed comparison between XPINNs and sE-PINN.
Figure~\ref{fig:xpinssandothers}  shows that XPINNs exhibits two different solutions on the interval node $t=0.5$, while the sE-PINN method strictly maintains continuity and smoothness. 
Table~\ref{HEXPINNSDEBIJIAO} shows that, compared to XPINNs, sE-PINN achieves better accuracy. Furthermore, using a single processor, sE-PINN requires less training time than XPINNs for this PDE.}

\begin{table}[h]
\color{black} 
\setlength{\abovecaptionskip}{0cm}
		\setlength{\belowcaptionskip}{0.2cm}
\caption{The performance of XPINNs and sE-PINN for solving Eq.~\eqref{llen-Cahn equation with strongly num:Allen-cahn}.}
\label{HEXPINNSDEBIJIAO}
\centering
\begin{adjustbox}{max width=\textwidth}
\begin{tabular}{c|c|c}
\hline 
Method & $\left\|\epsilon\right\|_2$  & Training Time  \\
\hline
XPINNs  & $  6.32 \times 10^{-3}$   & 15.68 min  \\ 
\hline
sE-PINN  & $  1.24 \times 10^{-3}$  & 13.87 min  \\  
\hline
\end{tabular}
\end{adjustbox}
\end{table}

\begin{figure}[!h]
    \centering    
\includegraphics[width=0.7\textwidth]{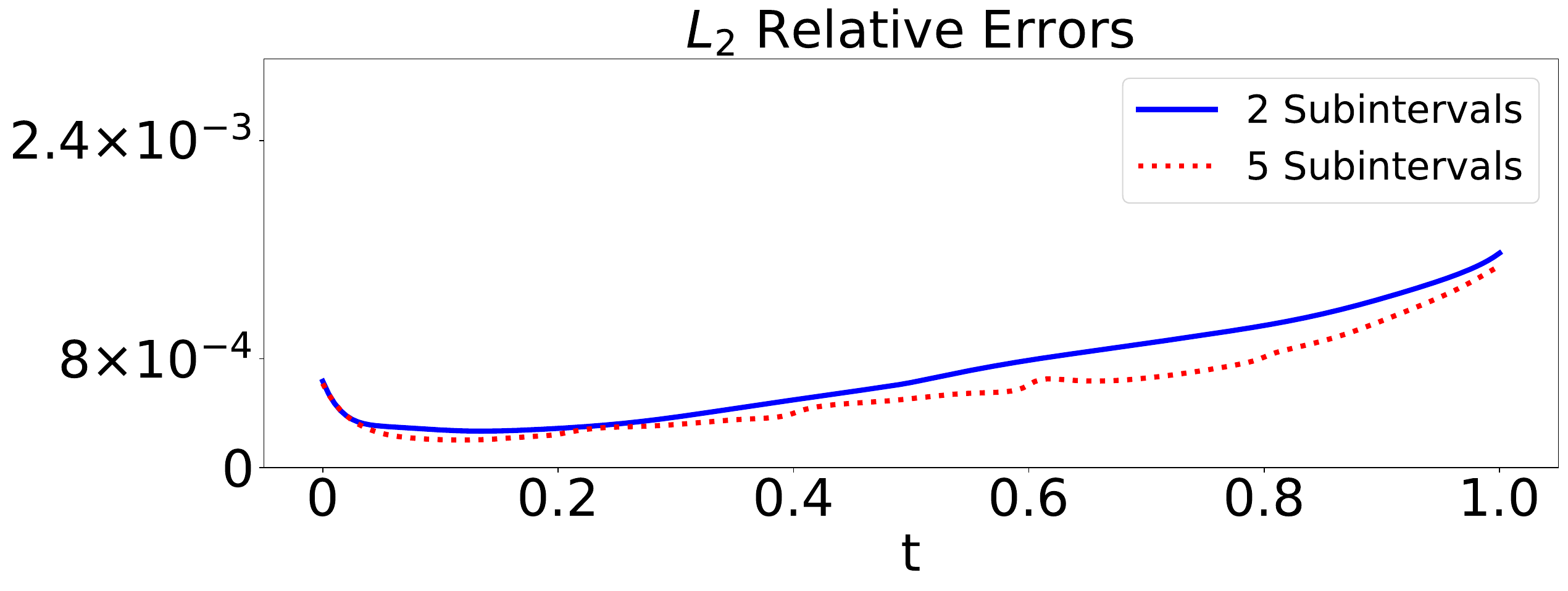}
    \caption{Variation curves of $L_2$ errors of aE-PINN methods with 2 and 5 intervals.} 
    \label{fig:acac}
\end{figure}

\begin{figure}[!h]
    \centering    
\includegraphics[width=0.9\textwidth]{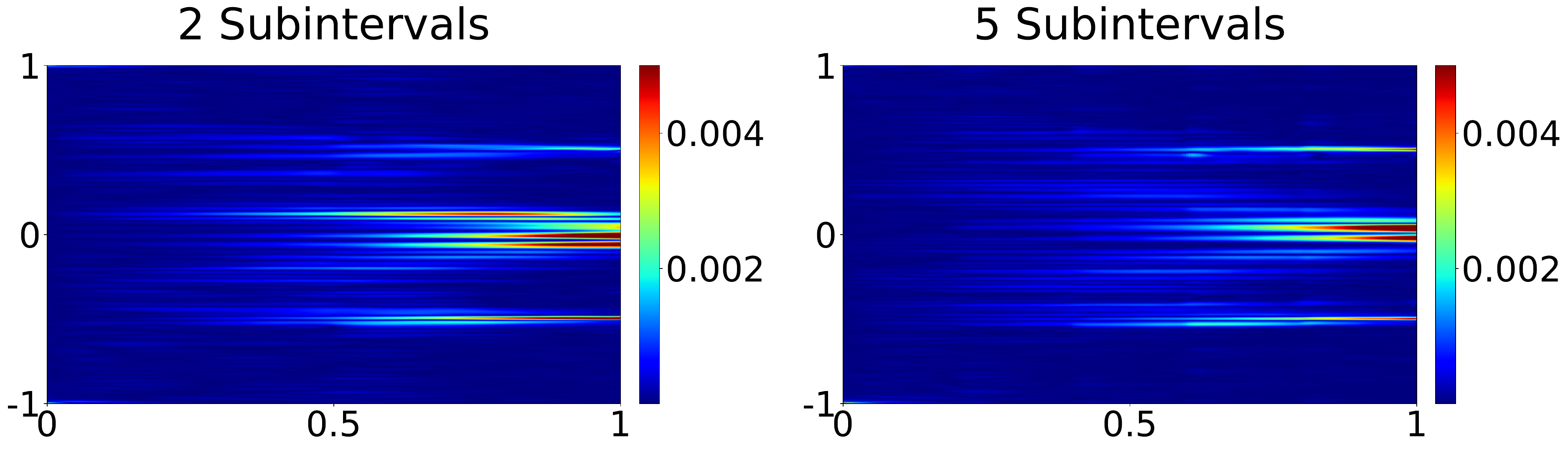}
    \caption{Point-wise errors of aE-PINN methods with 2 and 5 intervals for Eq.~\eqref{llen-Cahn equation with strongly num:Allen-cahn}.}
    \label{fig:acacac}
\end{figure}

In addition to the training results of 2 subintervals, we also list the result of 5 subintervals in Table~\ref{table Allen-Cahn equation with strongly nonlinear source term}.
Figure~\ref{fig:acac} shows the variation curves of $L_2$ errors of aE-PINN methods with 2 and 5 intervals, and
Figure~\ref{fig:acacac} shows their point-wise errors.
The training accuracy of 5 subintervals is slightly higher than that of 2 subintervals.

\subsection{Convection equation with large convection coefficient}\label{sec:Convection equation}

The convection equation~\eqref{consider Convection and convection equation} is commonly used to model some transport phenomena~\cite{krishnapriyan2021characterizing,haitsiukevich2023improved}. For large $\beta$, e.g., $\beta=40$, it is difficult to provide accurate PINN solution for the large time domain~\cite{PENWARDEN2023112464,krishnapriyan2021characterizing}.  
As discussed in Sec.~\ref{sec:Analysis of the extrapolation of PINNs}, this model has a weak extrapolation capability, so we use wE-PINN and aE-PINN to solve it to examine the performance of the two approaches.

We set 4 hidden layers for the E-DNN architecture, each layer with 100 neurons, the number of initial condition points $N_0=1200$, the number of  boundary condition points $N_b=1200$, and the number of residual points $N_r=10000$.

\begin{table}[h]
\color{black}
\setlength{\abovecaptionskip}{0cm}
		\setlength{\belowcaptionskip}{0.2cm}
\caption{The performance of different PINN versions for Eq.~\eqref{consider Convection and convection equation} with $\beta=40$.}
\label{Numerical Examples consider Convection and convection equation}
\centering
\begin{adjustbox}{max width=\textwidth}
\begin{tabular}{c|c|c}
\hline 
Method&  $\left\|\epsilon\right\|_2$  & Training Time \\
\hline
Curriculum Training~\cite{krishnapriyan2021characterizing}& $2.64 \times 10^{-3}$ &  78.48 min\\ 
\hline
wE-PINN (this work)&$4.45\times 10^{-3}$  &  64.57 min\\
\hline
aE-PINN (this work)&$5.65\times 10^{-4}$  &  68.05 min\\
\hline
\end{tabular}
\end{adjustbox}
\end{table}

Table~\ref{Numerical Examples consider Convection and convection equation} lists the results of different PINN versions. As 
\textcolor{black}{we can see}, our approaches have a better performance for this PDE as well.

\begin{figure}[ht]
    \centering    
    \includegraphics[width=0.99\textwidth]{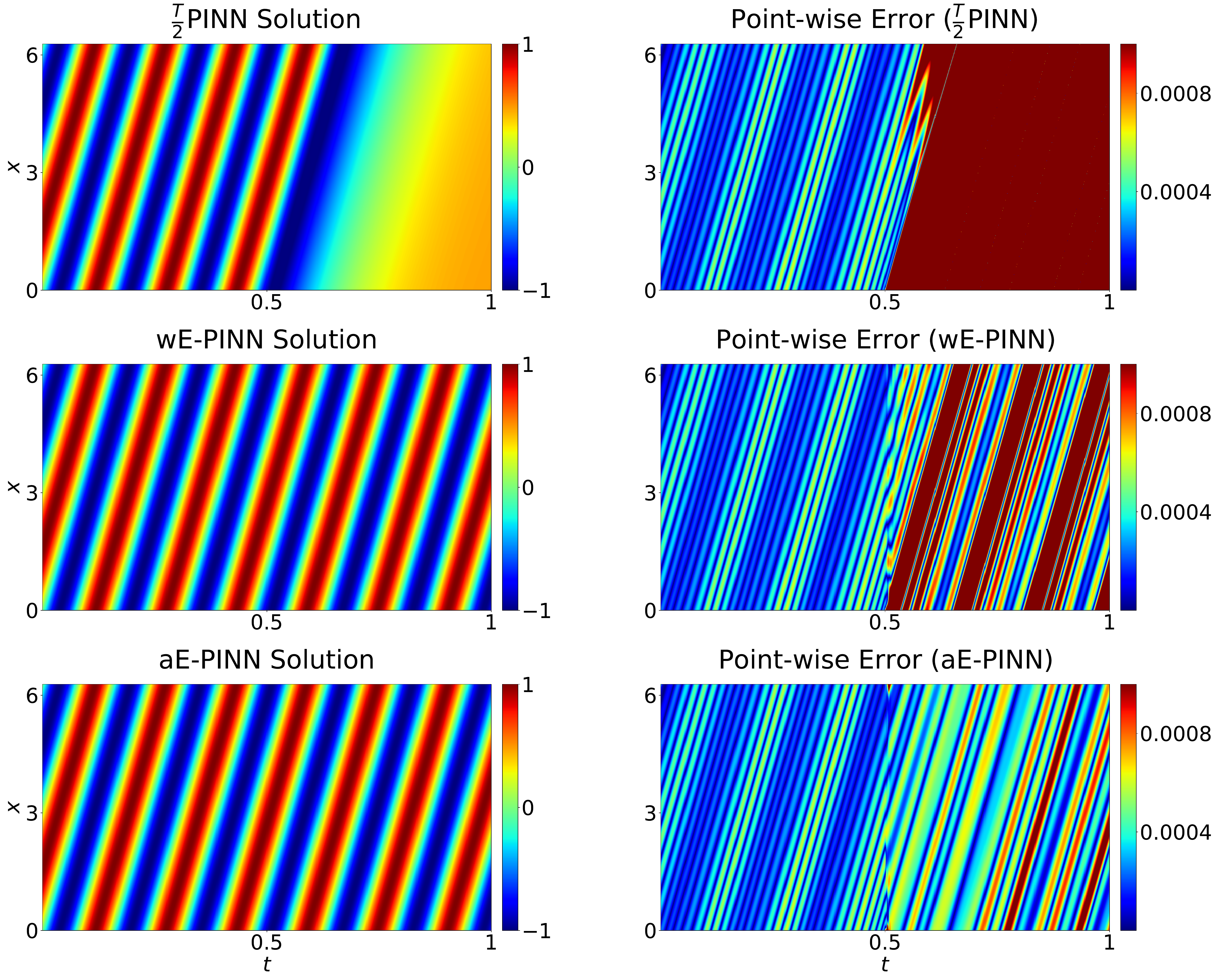}
    \caption{Numerical solutions and point-wise errors of different PINN versions for Eq.~\eqref{consider Convection and convection equation} with $\beta=40$.}
    \label{fig:Convection equation conventional pinn and s-ednn}
\end{figure}

Figure~\ref{fig:Convection equation conventional pinn and s-ednn} shows the images of solutions and point-wise errors for different PINN versions. 
If we focus only on the subinterval $[0,0.5]$, we can see that wE-PINN and aE-PINN have the same point-wise errors as that of the $\frac{T}{2}$PINN, which once again demonstrates that the E-DNN architecture successfully reproduces the solution of the previous subinterval.

The model has a weak extrapolation capability, so we have to use a rapidly varying function as the extrapolation control function to obtain an accurate solution.

\begin{figure}[!ht]
    \centering    
    \includegraphics[width=0.7\textwidth]{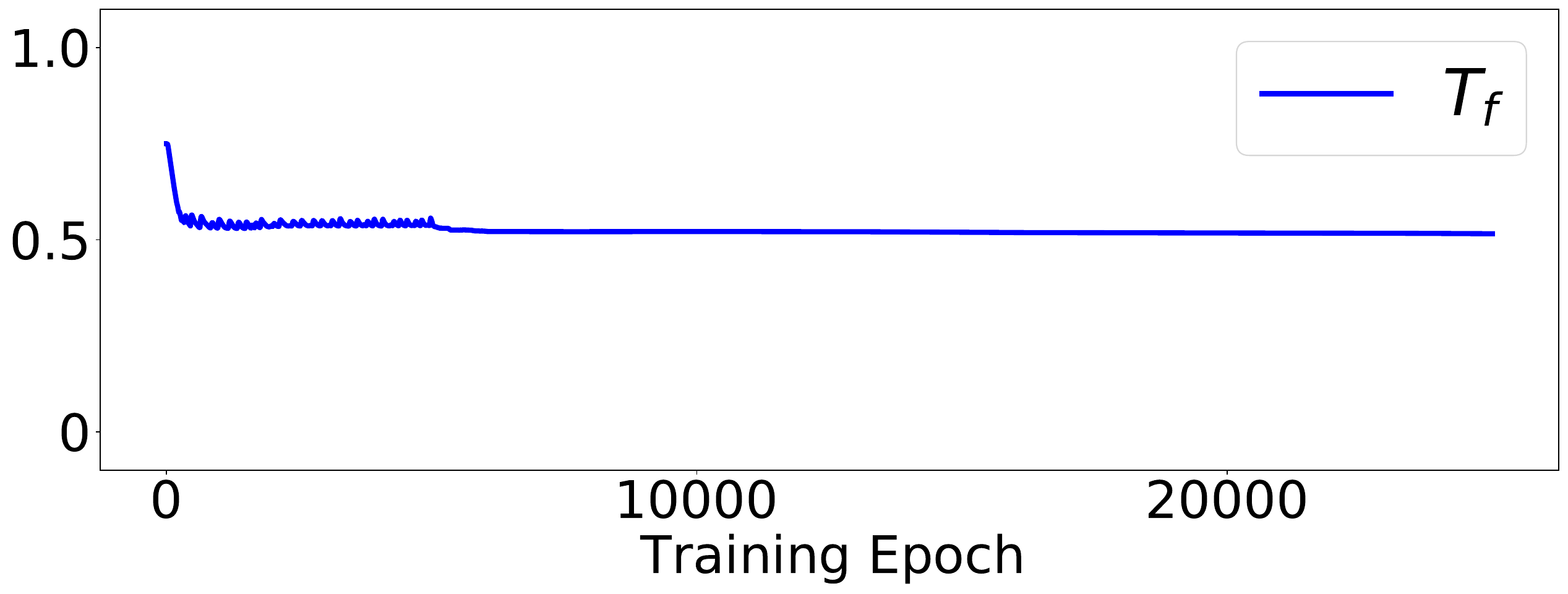}
    \caption{\textcolor{black}{Values of the trainable parameter $T_{f}$ in aE-PINN method for Eq.~\eqref{consider Convection and convection equation} with $\beta=40$.}}
    \label{The learnable parameters consider Convection and Convection convection equation}
\end{figure}

Figure~\ref{The learnable parameters consider Convection and Convection convection equation} shows the \textcolor{black}{values} of the trainable parameter $T_f$ during the optimization process when the adaptive control function $\mathcal{F}_a(t)$ is chosen. We can see that it is eventually optimized from the initial value of 0.75 to a value close to 0.5, which once again demonstrates the ability of the adaptive approach, it correctly identifies this model as one with weak extrapolation capabilities.

\begin{figure}[!ht]
    \centering    \includegraphics[width=0.95\textwidth]{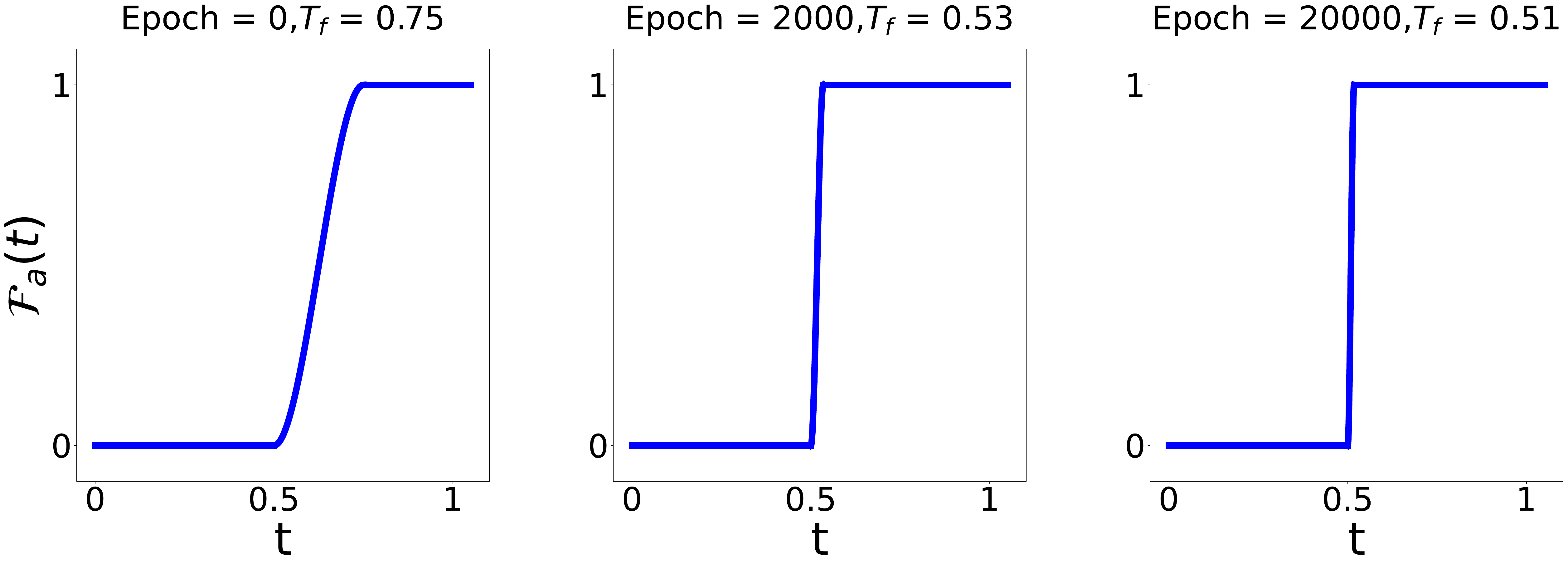}
    \caption{Optimization process of adaptive control function $\mathcal{F}_a(t)$ for Eq.~\eqref{consider Convection and convection equation} with $\beta=40$.}
    \label{duiliu2The learnable  parameters consider kdvkdvkdv equation}
\end{figure}

Figure~\ref{duiliu2The learnable  parameters consider kdvkdvkdv equation} shows the optimization process of the adaptive control function  $\mathcal{F}_a(t)$, it gradually becomes steeper from its initial smoother shape to fit the weak extrapolation property of this model.

\subsection{Korteweg-de Vries equation with high-order term}\label{sec:Korteweg-de vries equation}

This example is taken from the Ref.~\cite{PENWARDEN2023112464}, which describes the evolution of long one-dimensional waves in 
many physical settings~\cite{RAISSI1}.
The Korteweg-de Vries equation is formulated as follows,
\begin{equation}\label{num: let us consideKorteweg-de vries}
\begin{cases}
u_t+uu_{x}+0.0025u_{xxx}=0, &(x,t) \in (-1,1) \times  (0,1],\\
u(x,0)=\cos{\pi x},& x\in  (-1,1),   \\
u(-1,t)=u(1,t), u_x(-1,t)=u_x(1,t), & t\in  (0,1].\\
\end{cases}
\end{equation}


The E-DNN architecture used in this example contains 3 hidden layers, each layer with 30 neurons. We set $N_0=400$, $N_b=800$  and $N_r=8000$.

\begin{table}[h]
\color{black} 
\setlength{\abovecaptionskip}{0cm}
		\setlength{\belowcaptionskip}{0.2cm}
\caption{The performance of different PINN versions for Eq.~\eqref{num: let us consideKorteweg-de vries}.}
\label{Numerical Examples Korteweg-de Vries equation with high-order term}
\centering
\begin{adjustbox}{max width=\textwidth}
\begin{tabular}{c|c|c}
\hline 
Method&  $\left\|\epsilon\right\|_2$ & Training Time \\
\hline
s-d PINN~\cite{PENWARDEN2023112464}& $3.47\times 10^{-2}$ &310.93 min\\ 
\hline
aE-PINN (this work)& $6.20\times 10^{-3}$ & 23.15 min\\ 
\hline
\end{tabular}
\end{adjustbox}
\end{table}

Since we have no prior knowledge of the extrapolation capability for this PDE, we choose the adaptive interpolation-driven architecture for it. 
Table~\ref{Numerical Examples Korteweg-de Vries equation with high-order term} lists the results of different PINN versions.

\begin{figure}[!ht]
    \centering    
    \includegraphics[width=0.7\textwidth]{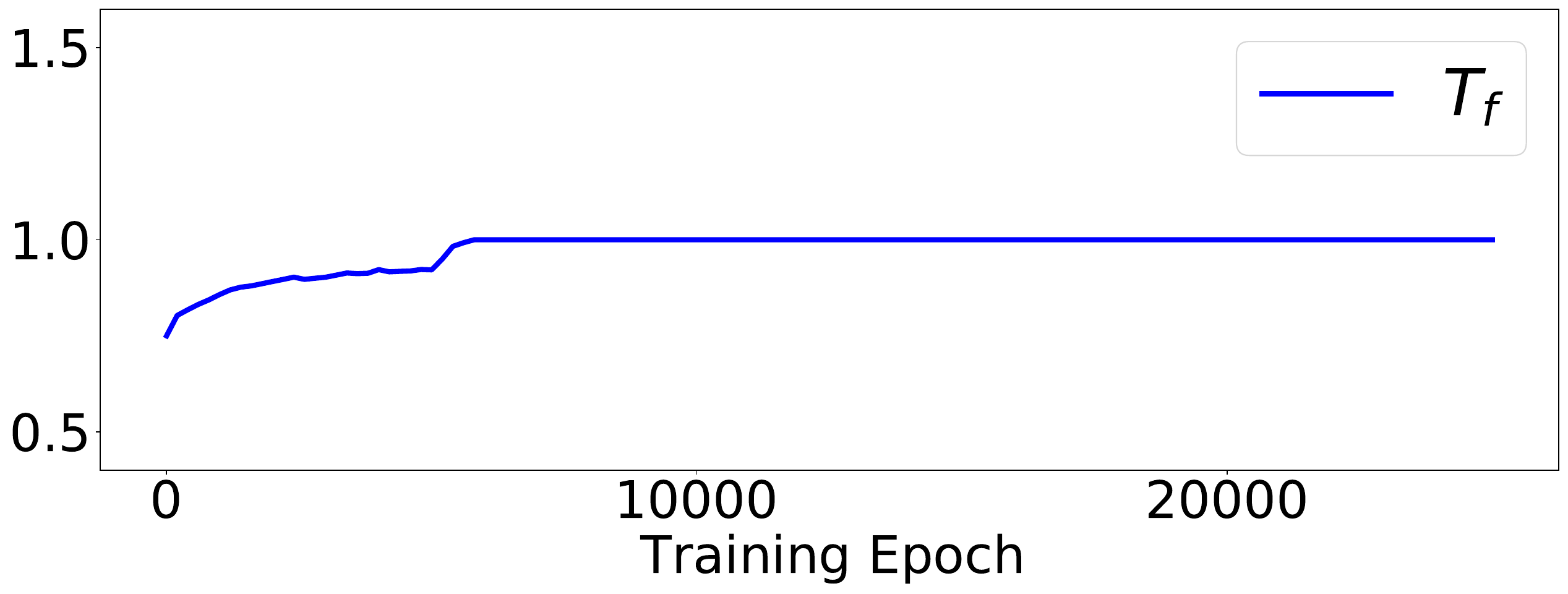}
    \caption{\textcolor{black}{Values of the trainable parameter $T_{f}$ in aE-PINN method for Eq.~\eqref{num: let us consideKorteweg-de vries}}.}
    \label{The KDVVVlearnable parameters consider kdvkdvkdv equation}
\end{figure}

\begin{figure}[!ht]
    \centering    
\includegraphics[width=0.9\textwidth]{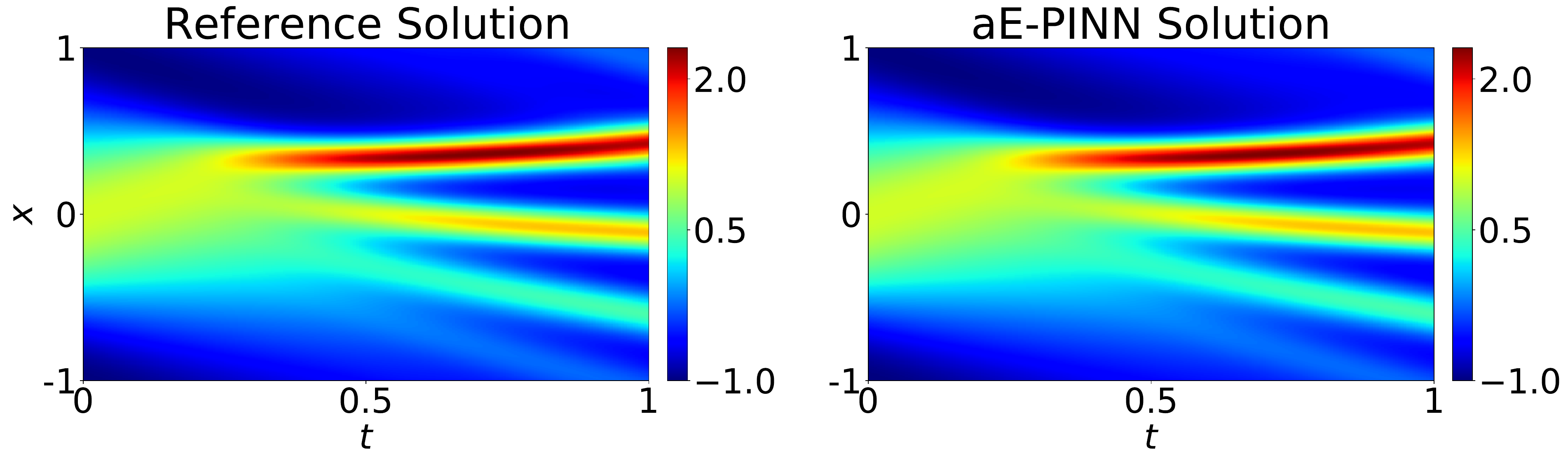}
    \caption{The reference solution and aE-PINN solution of Eq.~\eqref{num: let us consideKorteweg-de vries}.}
    \label{fig:acacacThe sE-PINN and aE-PINN method gives the best prediction among them}
\end{figure}

Figure~\ref{The KDVVVlearnable parameters consider kdvkdvkdv equation} shows the \textcolor{black}{values} of the trainable parameter $T_f$ during the optimization process. It is eventually optimized from 0.75 to 1.0, indicating that this model has a strong extrapolation capability.
Figure~\ref{fig:acacacThe sE-PINN and aE-PINN method gives the best prediction among them} shows the images of the reference solution and aE-PINN solution.

\section{Conclusion}
\textcolor{black}{To address the limitations of sequential learning strategies,}
we systematically \textcolor{black}{analyzed} the extrapolation property of the PINN solution, then motivated by the extrapolation capability, we propose a new neural network architecture called E-DNN. The new PINN approaches derived from the E-DNN architecture can perfectly address the aforementioned problems. Moreover, our approaches are relatively computationally efficient and have inherent scalability. With the proposed approaches, using only a single neural network, we can divide the large time domain into several subintervals and then obtain the local PINN solution of each of the subintervals one by one, and the training process is independent among them; especially, the overall solution of the whole time domain combined by these local PINN solutions is seamless, it strictly preserves the continuity and smoothness at the interval nodes.
In addition, our approach has good scalability and is naturally suited for high-dimensional problems. However, it also has some limitations. One of the potential problems is that in the process of one-by-one training, the error of the previous time interval could be \textcolor{black}{propagated into} the next time interval, resulting in a certain amount of error accumulation. How to reduce or eliminate the possible error accumulation is our future work.


\section*{Acknowledgements}
\textcolor{black}{
The authors sincerely thank the reviewers for their numerous insightful comments and suggestions, which have greatly improved the presentation of the paper.}
The work of Y. Yao is supported by the National Natural Science Foundation of China (No. 12271055). The work of Z. Gao is supported by the National Key R$\&$D Program of China (No. 2022YFA1004500), the National Natural Science Foundation of China (No. 12171048) and the Funding of the National Key Laboratory of Computational Physics. 
\appendix
\section{Testing results of the effect of training data and number of neurons on the extrapolation capability of PINNs.}\label{appendixaaaaaaaa}
\begin{figure}[h]
    \centering    
        \includegraphics[width=0.84\textwidth]{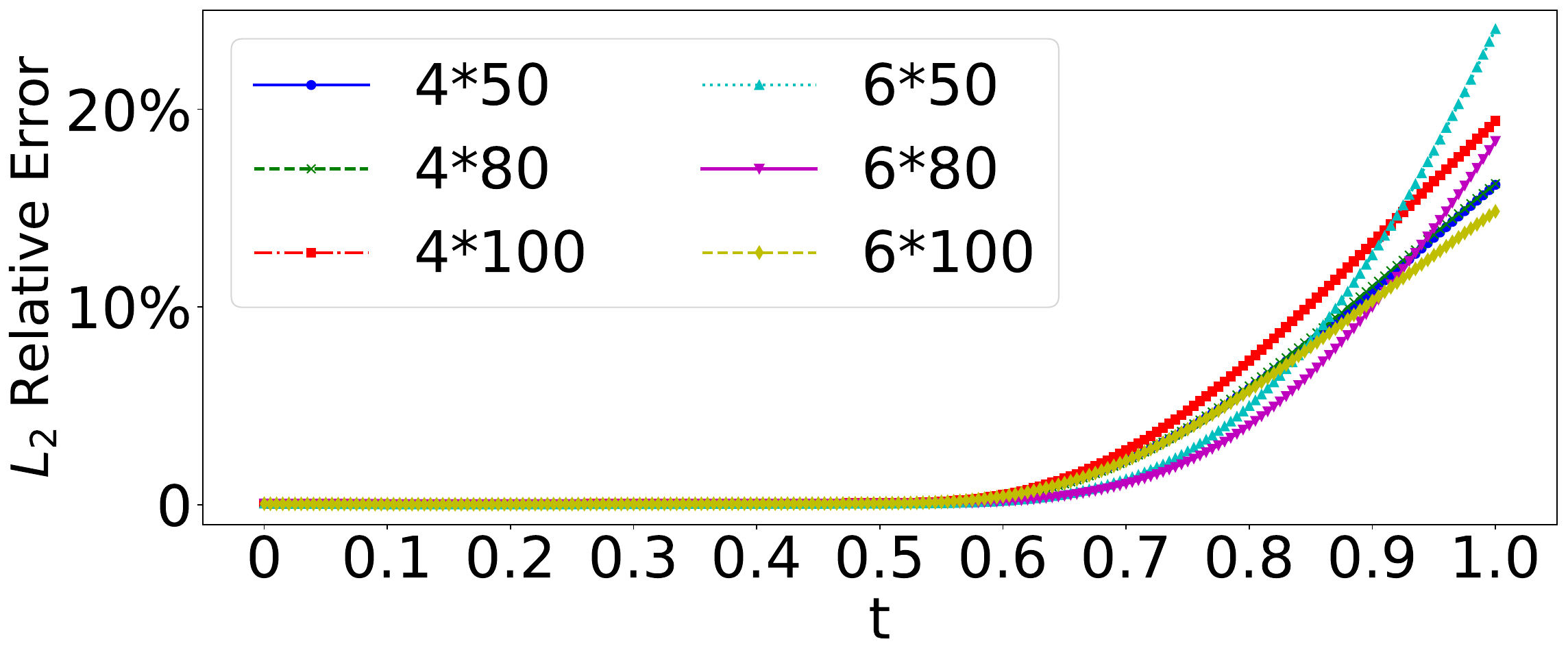}
    \caption{$L_2$ errors of $\frac{T}{2}$PINN with different hidden layers for Eq.~\eqref{llen-Cahn equation with strongly num:Allen-cahn}.} 
    \label{fig:comssdasdadasdasssv}
\end{figure}

\begin{figure}[h]
    \centering    
        \includegraphics[width=0.84\textwidth]{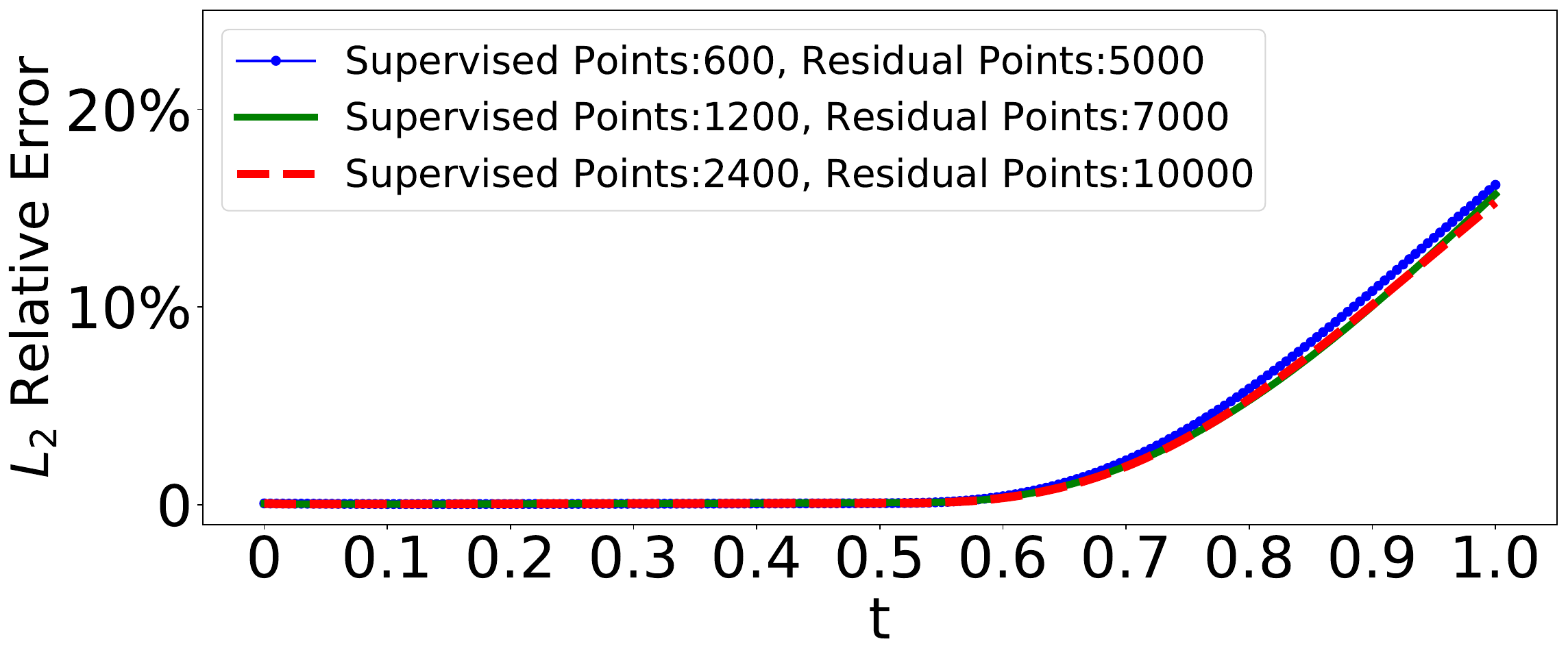}
    \caption{$L_2$ errors of $\frac{T}{2}$PINN  with different number of sampling points for Eq.~\eqref{llen-Cahn equation with strongly num:Allen-cahn}.} 
    \label{fig:comsfsdgsdagssv}
\end{figure}

\begin{figure}[!h]
    \centering    
        \includegraphics[width=0.84\textwidth]{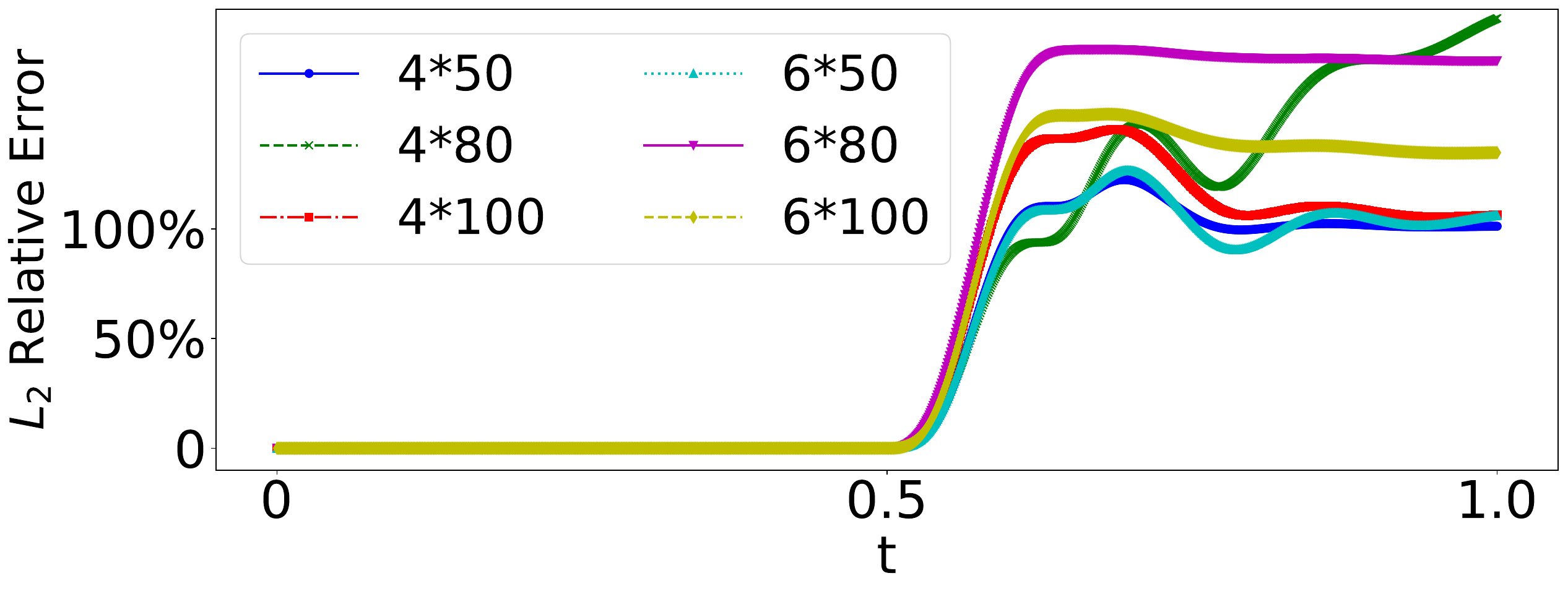}
    \caption{$L_2$ errors of $\frac{T}{2}$PINN  with different hidden layers  for Eq.~\eqref{consider Convection and convection equation} with $\beta=40$.}
    \label{fig:duiliucomssdasdadasdasssv}
\end{figure}

\begin{figure}[h]
    \centering    
    \includegraphics[width=0.84\textwidth]{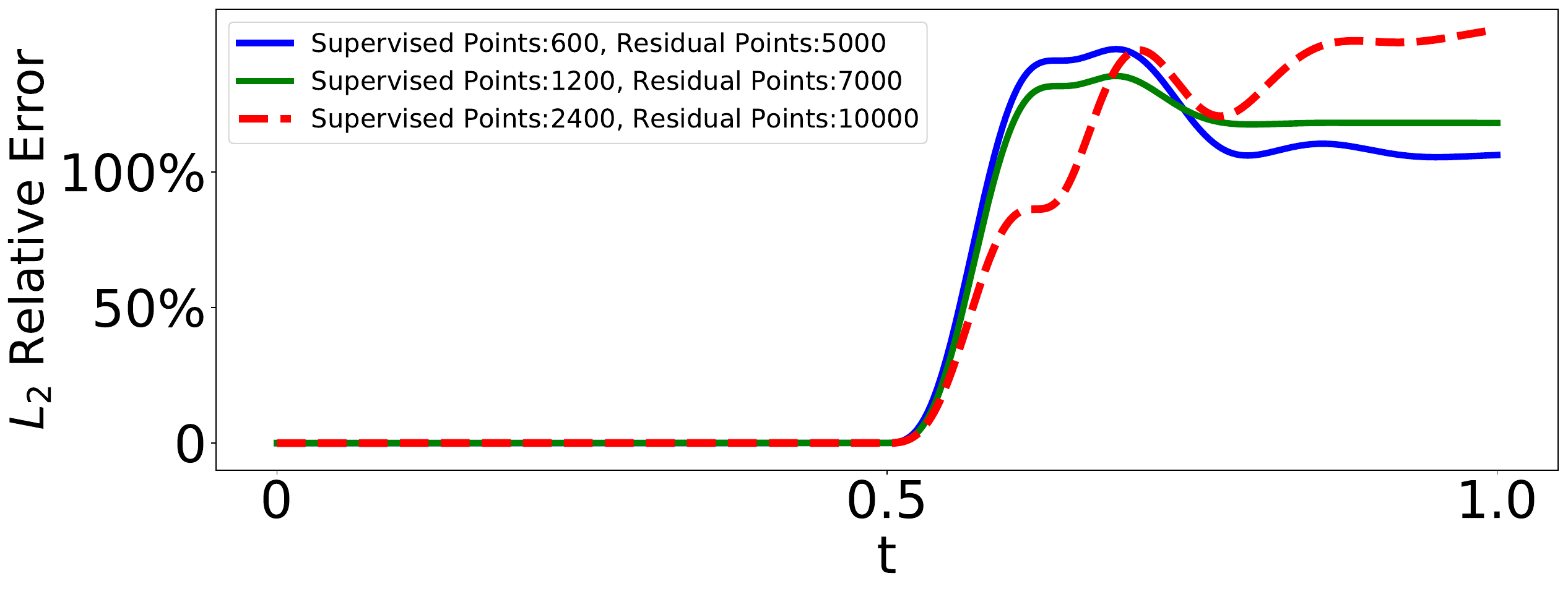}
    \caption{$L_2$ errors of $\frac{T}{2}$PINN  with different number of sampling points  for Eq.~\eqref{consider Convection and convection equation} with $\beta=40$.}
    \label{fig:duiliucomsfsdgsdagssv}
\end{figure}

For Eq.~\eqref{llen-Cahn equation with strongly num:Allen-cahn},  Figure~\ref{fig:comssdasdadasdasssv} shows the extrapolation performance of $\frac{T}{2}$PINN obtained with different network sizes for a fixed number of training points. Figure~\ref{fig:comsfsdgsdagssv} gives the results of different numbers of training points for a fixed neural network size.
Figures~\ref{fig:duiliucomssdasdadasdasssv} and \ref{fig:duiliucomsfsdgsdagssv} \textcolor{black}{display} the similar results for Eq.~\eqref{consider Convection and convection equation} with $\beta=40$.

\section{Various types of extrapolation control functions}\label{appendixBBBBBBBB}
\textcolor{black}{
Table~\ref{Different control function} lists 5 control functions for $\mathcal{F}_s(t)$ that satisfy the conditions described in Section~\ref{sec:An strong}. 
Table~\ref{duozhongkongzhihanshu} shows the accuracy, training time and number of iteration steps of the sE-PINN method for solving Eq.~\ref{llen-Cahn equation with strongly num:Allen-cahn}. We can see that the choice of the control function does not significantly affect the accuracy and efficiency of the sE-PINN method.}

\begin{table}[!h]
\setlength{\abovecaptionskip}{0cm}
		\setlength{\belowcaptionskip}{0.2cm}
\centering
\caption{Five control functions $\mathcal{F}_s(t)$ for sE-PINN.} 
\label{Different control function}
\begin{adjustbox}{max width=\textwidth}
\begin{tabular}{c|c|c|c}
\hline
\diagbox{$\mathcal{F}_{s}(t)$}{$t$}& $[0,0.5)$  & $[0.5,1)$ & $[1,+\infty)$ \\
\hline
$\mathcal{F}_{s1}(t)$ & 0 & $ \cos^2(\pi t)$ & 1 \\
\hline
$\mathcal{F}_{s2}(t)$ &  0 & $- 16 t^{3} + 36 t^{2} - 24t + 5$  & 1 \\
\hline
$\mathcal{F}_{s3}(t)$ & 0 & $16 t^{4} - 64 t^{3} + 88 t^{2} - 48 t + 9$ & 1 \\
\hline
$\mathcal{F}_{s4}(t)$ & 0 & $e^{4 t^{2} - 6 t + 2} - 16 t^{3} + 32 t^{2} - 18 t  + 2$ & 1 \\
\hline
$\mathcal{F}_{s5}(t)$ & 0 & $e^{16 t^{4} - 48 t^{3} + 52 t^{2} - 24 t + 4} + \cos^{2}{(\pi t )} - 1$ & 1 \\
\hline
\end{tabular}
\end{adjustbox}
\end{table}

\begin{table}[!h]
\setlength{\abovecaptionskip}{0cm}
		\setlength{\belowcaptionskip}{0.2cm}
\caption{Results of sE-PINN with various control function 
$\mathcal{F}_s(t)$ for solving Eq.~\ref{llen-Cahn equation with strongly num:Allen-cahn}.}
\label{duozhongkongzhihanshu}
\centering
\begin{adjustbox}{max width=\textwidth}
\begin{tabular}{c|c|c|c}
\hline 
$\mathcal{F}_s(t)$ & $\left\|\epsilon\right\|_2$  & Training Time & Iterations\\
\hline
$\mathcal{F}_{s1}(t)$ & $  1.25 \times 10^{-3}$   & 13.92 min & 24330 \\ 
\hline
$\mathcal{F}_{s2}(t)$  & $  1.24 \times 10^{-3}$  & 13.87 min & 24358\\ 
\hline
$\mathcal{F}_{s3}(t)$  & $  1.21 \times 10^{-3}$   & 13.83 min & 24219 \\ 
\hline
$\mathcal{F}_{s4}(t)$  & $  1.24 \times 10^{-3}$  & 13.90 min  & 24449\\ 
\hline
$\mathcal{F}_{s5}(t)$   & $  1.23 \times 10^{-3}$   & 14.03 min &  24480\\ 
\hline
\end{tabular}
\end{adjustbox}
\end{table}


\textcolor{black}{
Table~\ref{wwwDifferent control function} lists 5 control functions for $\mathcal{F}_w(t)$, 
Table~\ref{wwwduozhongkongzhihanshu} lists the results of the wE-PINN method with different control functions for solving Eq.~\eqref{consider Convection and convection equation} with $\beta=40$. 
They all reach the maximum number of iterations set for the optimizer.}

\begin{table}[!h]
\setlength{\abovecaptionskip}{0cm}
		\setlength{\belowcaptionskip}{0.2cm}
\centering
\caption{Five control functions $\mathcal{F}_w(t)$ for wE-PINN.} 
\label{wwwDifferent control function}
\begin{adjustbox}{max width=\textwidth}
\begin{tabular}{c|c|c|c}
\hline
\diagbox{$\mathcal{F}_{w}(t)$}{$t$}& $[0,0.5)$  & $[0.5,0.6)$ & $[0.6,+\infty)$ \\
\hline
$\mathcal{F}_{w1}(t)$ & 0 & $ \cos^2(5\pi t)$ & 1 \\
\hline
$\mathcal{F}_{w2}(t)$ &  0 & $- 2000 t^{3} + 3000 t^{2} - 1800t + 325$   & 1 \\
\hline
$\mathcal{F}_{w3}(t)$ & 0 & $10000 t^{4} - 24000 t^{3} + 21400 t^{2} - 8400 t + 1225$ & 1 \\
\hline
$\mathcal{F}_{w4}(t)$ & 0 & $e^{100 t^{2} - 110 t + 30} - 2000 t^{3} + 3200 t^{2} - 1690 t  + 294$ & 1 \\
\hline
$\mathcal{F}_{w5}(t)$ & 0 & $e^{10000 t^{4} - 22000 t^{3} + 18100 t^{2} - 6600 t + 900} + \cos^{2}{(5\pi t )} - 1$ & 1 \\
\hline
\end{tabular}
\end{adjustbox}
\end{table}

\begin{table}[!h]
\setlength{\abovecaptionskip}{0cm}
		\setlength{\belowcaptionskip}{0.2cm}
\caption{Results of wE-PINN with various control functions 
$\mathcal{F}_w(t)$ for solving Eq.~\eqref{consider Convection and convection equation} with $\beta=40$.}
\label{wwwduozhongkongzhihanshu}
\centering
\begin{adjustbox}{max width=\textwidth}
\begin{tabular}{c|c|c|c}
\hline 
$\mathcal{F}_w(t)$ & $\left\|\epsilon\right\|_2$  & Training Time & Iterations\\
\hline
$\mathcal{F}_{w1}(t)$ & $  5.52 \times 10^{-3}$   & 63.10 min & 25000 \\ 
\hline
$\mathcal{F}_{w2}(t)$  & $  4.45 \times 10^{-3}$  & 64.57 min & 25000\\ 
\hline
$\mathcal{F}_{w3}(t)$  & $  4.53 \times 10^{-3}$   & 63.87 min & 25000 \\ 
\hline
$\mathcal{F}_{w4}(t)$  & $  3.53 \times 10^{-3}$  & 64.13 min  & 25000\\ 
\hline
$\mathcal{F}_{w5}(t)$   & $  4.41 \times 10^{-3}$   & 64.45 min &  25000\\ 
\hline
\end{tabular}
\end{adjustbox}
\end{table}

\end{document}